\documentclass[11pt]{article}
\usepackage{amsmath,amssymb}
\usepackage[dvipdfmx]{graphicx}
\begin{document}

\title{On the triplicity among infinite products, infinite series, and continued fractions;\\ and its 
applications to divergent series}
\author{Kiyoshi \textsc{Sogo}
\thanks{EMail: sogo@icfd.co.jp}
}
\date{}

\maketitle

\begin{center}
Institute of Computational Fluid Dynamics, 1-16-5, Haramachi, Meguroku, 
Tokyo, 152-0011, Japan
\end{center}
\abstract{
Many identities written by $P=S=C$ are obtained, where $P$ infinite products, $S$ infinite series, 
and $C$ continued fractions. Such equality is called {\it triplicity}, and it can be used to compute 
the values of infinite series. It is applied even to obtain sums of divergent series. 
Many examples of such infinite series are shown, including $1-2+2^3-2^6+\cdots$, 
which is in Entry 7 of Gauss's diary and its value $0.4275251302\cdots$ is here obtained.
}
%
%
\tableofcontents

\section{Introduction}
\setcounter{equation}{0}

\subsection{Concepts of triplicity and pivot variable}

The {\it triplicity}, or triad, introduced here is well exemplified by the famous Cauchy's identity \cite{Andrews} among 
infinite product ($P$), infinite series ($S$), and infinite continued fraction ($C$), in the form 
of identity $P=S=C$, that is
\begin{align}
\prod_{n=0}^\infty \frac{(1-\alpha xq^n)}{(1-xq^n)}&=\sum_{n=0}^\infty \frac{(\alpha; q)_n}{(q; q)_n}\ x^n
\nonumber \\
&=\frac{e_0}{1-}\ \frac{e_1x}{1-}\ \frac{e_2x}{1-}\ \frac{e_3x}{1-}\ \cdots,
\label{Cauchy}
\end{align}
where we require $|q|<1,\ |x|<1$ for convergence, and we set $q$-Pochhammer symbol
\begin{align}
\begin{split}
&(\alpha; q)_0=1	,\\
&(\alpha; q)_n=(1-\alpha)(1-\alpha q)\cdots (1-\alpha q^{n-1}),\quad (n \geq 1).
\end{split}
\label{qPochhammer}
\end{align}
The coefficients $e_0=1,\ e_1,\ e_2,\ e_3,\ \cdots$, are 
\begin{align}
&e_1=\frac{1-\alpha}{1-q},\ e_2=\frac{\alpha-q}{(1+q)(1-q)},\ e_3=\frac{q (1-\alpha q)}{(1+q)(1-q^3)},\ \cdots,
\end{align}  
while general expressions are given for positive integer $m\geq 1$ by
\begin{align}
e_{2m}=\frac{q^{m-1}(\alpha-q^m)}{(1+q^m)(1-q^{2m-1})},\quad 
e_{2m+1}=\frac{q^{m}(1-\alpha q^m)}{(1+q^m)(1-q^{2m+1})},
\end{align}
which are derived by using formulas found by Muir \cite{Muir} and Rogers \cite{Rogers}, which will be  
given in the next section by \eqref{MuirF} and \eqref{MuirRogersF}. 
For the above case \eqref{MuirRogersF} is used.

Let us call the variable $x$ in \eqref{Cauchy} the {\it pivot} variable, because it works as a generating 
parameter.
Some triplicity identities however do not have such variable. 
Typical example is the well-known Rogers-Ramanujan identity
\begin{align}
\prod_{n=0}^\infty \frac{(1-q^{5n+1})(1-q^{5n+4})}{(1-q^{5n+2})(1-q^{5n+3})}&=
\frac{1+\sum_{n=1}^\infty \frac{q^{n(n+1)}}{(q; q)_n}}
{1+\sum_{n=1}^\infty \frac{q^{n^2}}{(q; q)_n}} \nonumber \\
&=\frac{1}{1+}\ \frac{q}{1+}\ \frac{q^2}{1+}\ \frac{q^3}{1+}\ \cdots,
\label{RRIdentity}
\end{align}
which has no $x$ variable. Here it should be noted that we include the {\it quotient} of two 
infinite series also in $S$ group. 
For this example we can actually introduce variable $x$ {\it partially} in $S=C$ part, such as
\begin{align}
\frac{1+\sum_{n=1}^\infty \frac{q^{n(n+1)}}{(q; q)_n} x^n}
{1+\sum_{n=1}^\infty \frac{q^{n^2}}{(q; q)_n} x^n}
=\frac{1}{1+}\ \frac{q x}{1+}\ \frac{q^2 x}{1+}\ \frac{q^3 x}{1+}\ \cdots.
\end{align}
This equality can also be derived by Muir formula \eqref{MuirF} in the next section. 
However it is not known unfortunately how to include variable $x$ in the infinite product part $P$. 
Therefore such case might be called as {\it weak-triplicity}, and the case like \eqref{Cauchy} 
that all $P,\ S,\ C$ contain pivot variable, as {\it strong-triplicity} for the sake of distinction. 

\subsection{Triplicity can be applied to divergent series}

According to Slater (in the first page of \cite{Slater}), Wallis in 1655 considered 
a divergent series
\begin{align}
\sum_{n=0}^\infty (-1)^n n!=1-1!+2!-3!+\cdots,
\label{WallisSum}
\end{align}
which is a special case ($a=1,\ x=1$) of
\begin{align}
\sum_{n=0}^\infty (-1)^n (a)_n\ x^n,
\label{WF}
\end{align}
where $(a)_n$ is the Pochhammer symbol
\begin{align}
(a)_0=1,\quad (a)_n=a(a+1)\cdots(a+n-1),\quad (n\geq 1).
\label{Pochhammer}
\end{align}
The alternating power series \eqref{WF} is generically still divergent except $a$ is a negative integer, 
since it has zero radius of convergence.  
In fact it is an example of asymptotic power series.

If $P$ part is replaced by definite integral, triplicity for the Wallis problem is given  
with the pivot variable $x>0$ by
\begin{align}
\frac{1}{\Gamma(a)}\ \int_0^\infty \frac{s^{a-1}\ e^{-s}}{1+xs}\ ds
&=\sum_{n=0}^\infty (-1)^n (a)_n x^n
\nonumber \\
&=\frac{e_0}{1-}\ \frac{e_1x}{1-}\ \frac{e_2 x}{1-}\ \frac{e_3 x}{1-}\ \cdots,
\label{WallisTriplicity}
\end{align}
where the coefficients in the second line are 
\begin{align}
e_0=1,\ e_1=-a,\ e_2=-1,\ e_3=-(a+1),\ \cdots,
\end{align}
whose general expressions are given by
\begin{align}
e_{2m}=-m,\quad e_{2m+1}=-(a+m),\qquad (m \geq 1),
\end{align}
which are derived by using Muir-Rogers formula \eqref{MuirRogersF} below. 

For the sake of completeness, let us give a {\it proof} of the above equality. 
The integral in the left hand side is rewritten as follows,
\begin{align}
&\frac{1}{\Gamma(a)}\int_0^\infty s^{a-1}\ e^{-s} \left( \int_0^\infty e^{-t(1+xs)} dt \right) ds
\nonumber \\
&=\frac{1}{\Gamma(a)} 
\sum_{n=0}^\infty (-1)^n \frac{x^n}{n!}\left(\int_0^\infty  s^{a+n-1} e^{-s} ds\int_0^\infty t^n e^{-t} dt \right) 
\nonumber \\
&=\sum_{n=0}^\infty (-1)^n x^n\cdot \frac{\Gamma(a+n)}{\Gamma(a)}
=\sum_{n=0}^\infty (-1)^n (a)_n x^n,
\nonumber
\end{align}
because $\Gamma(a+n)/\Gamma(a)=(a)_n$ due to the iterative property of Gamma function. 

It should be added one remark. Why the converging integral for $a>0$ in the left hand side changes into a divergent series? 
The reason is due to an exchange of integration and summation, which is not always allowed. 
In other word, here we {\it define} the right hand side series by the left hand side integral. 
Such procedure can also be said a kind of {\it analytical continuation}.

Therefore we can compute the Wallis sum by numerical integration
\begin{align}
\sum_{n=0}^\infty (-1)^n n!=\int_0^\infty \frac{e^{-s}}{1+s}\ ds=0.596347\cdots,
\end{align}
replacing the divergent series by the integral in triplicity relation. 
Methods of computing divergent series, including the above, is also discussed by Hardy \cite{HardyD} for example.

In a similar manner, we can obtain sums of divergent series by using triplicity relations. 
Early examples are given by Gauss in Entry 7 (1796) of his mathematical diary \cite{GaussD}, {\it without answers},
\begin{align}
&(1)\quad 1-2+8-64+\cdots=\frac{1}{1+}\ \frac{2}{1+}\ \frac{2}{1+}\ \frac{8}{1+}\ \frac{12}{1+}\ \cdots,
\\
&(2)\quad 1-1+1\cdot 3-1\cdot 3\cdot 7+\cdots
=\frac{1}{1+}\ \frac{1}{1+}\ \frac{2}{1+}\ \frac{6}{1+}\ \cdots,
\end{align}
which are the special case ($q=2$) of
\begin{align}
&(1)\quad \sum_{n=0}^\infty (-1)^n q^{n(n+1)/2}=\frac{1}{1+}\ \frac{q}{1+}\ \frac{q^2-q}{1+}\ \frac{q^3}{1+}\ \cdots,
\label{Gauss_1}
\\
&(2)\quad \sum_{n=0}^\infty (q; q)_n
=\frac{1}{1-}\ \frac{1-q}{1-}\ \frac{q-q^2}{1-}\ \frac{q-q^3}{1-}\ \cdots,
\label{Gauss_2}
\end{align}
both sides of which are divergent for $q>1$. \\

The composition of the present article is as follows. 
For the preparations we first discuss in \S 2 various methods making continued fractions from a given series. 
After such preparations we will discuss in \S 3 about various triplicity and apply them to obtain their sums. 
Problems computing sums of infinite series, including divergent cases, are the main subjects of \S 3, 
where above questions by Gauss are answered together with other problems. 


It should be remarked here that Ramanujan gave his method of computing divergent series in \S 6 of his Notebook 
\cite{Ramanujan}, based mainly on Euler-MacLaurin type formula. 
The procedure of the present article based on the triplicity among $P,\ S$ and $C$, is different from Ramanujan's 
method described there. 

The present idea introducing $p=q^{-1}$ for $q>1$, is inspired by the discussion of Slater \cite{Slater} 
on convergence conditions of $q$-hypergeometric series. 
And many examples in \S 2 and especially in \S 3.1, which are used for our discussion on divergent series, 
are borrowed from an article by Andrews \cite{Andrews2}, which surveys several topics in so-called {\it Lost Notebook} of Ramanujan. 


\section{Methods of making continued fractions}
\setcounter{equation}{0}

\subsection{Euler formula for continued fractions}

{\bf (Continued fractions)}\quad 
In his book {\it Introductio in Analysin Infinitorum} (Introduction to Analysis of the Infinite, 1748) \cite{Euler}, 
Euler gave his famous expression of infinite series by continued fraction such that
\begin{align}
c_0-c_1+c_2-c_3+\cdots
&=\frac{a_0}{b_0+}\ \frac{a_1}{b_1+}\ \frac{a_2}{b_2+}\ \frac{a_3}{b_3+}\ \cdots,
\label{EulerCF}
\end{align}
where he found 
\begin{align}
&a_0=c_0,\ b_0=1;\quad a_1=c_1,\ b_1=c_0-c_1;\nonumber \\
&a_2=c_0c_2,\ b_2=c_1-c_2;\quad
a_3=c_1c_3,\ b_3=c_2-c_3;\ \cdots,
\end{align}
which are given generally by
\begin{align}
a_n=c_{n-2}c_n,\ b_n=c_{n-1}-c_n,\qquad (n\geq 2).
\end{align}

Euler probably considered as follows to derive his formula. Firstly let us rewrite 
\begin{align}
c_0-c_1=\frac{c_0(c_0-c_1)}{c_0}=\frac{c_0}{\frac{c_0}{c_0-c_1}}=\frac{c_0}{1+\frac{c_1}{c_0-c_1}}.
\nonumber
\end{align}
Then replacing $c_1$ by $c_1-c_2$ and rewrite
\begin{align}
\frac{c_1}{c_0-c_1}\rightarrow \frac{c_1-c_2}{c_0-c_1+c_2}=\frac{c_1(c_1-c_2)}{c_1(c_0-c_1+c_2)}
=\frac{c_1}{c_0-c_1+\frac{c_0c_2}{c_1-c_2}},
\nonumber
\end{align}
where we used $c_1(c_0-c_1+c_2)=(c_0-c_1)(c_1-c_2)+c_0c_2$. Then replacing $c_2$ by $c_2-c_3$, and rewrite 
\begin{align}
\frac{c_0c_2}{c_1-c_2}\rightarrow \frac{c_0(c_2-c_3)}{c_1-c_2+c_3}=\frac{c_0c_2(c_2-c_3)}{c_2(c_1-c_2+c_3)}
=\frac{c_0c_2}{c_1-c_2+\frac{c_1c_3}{c_2-c_3}},
\nonumber
\end{align}
where we used $c_2(c_1-c_2+c_3)=(c_1-c_2)(c_2-c_3)+c_1c_3$. Such procedure can be repeated infinitely,  and 
we obtain Euler formula \eqref{EulerCF}.

We can make the correspondence between series and continued fraction in the unique form 
(one to one correspondence) by writing
\begin{align}
&c_0-c_1+c_2-c_3+\cdots = \frac{d_0}{1+}\ \frac{d_1}{1+}\ \frac{d_2}{1+}\ \frac{d_3}{1+}\ \cdots,
\end{align}
where we have defined
\begin{align}
&d_0=\frac{a_0}{b_0}=c_0,\quad d_1=\frac{a_1}{b_0b_1}=\frac{c_1}{c_0-c_1},\quad 
d_2=\frac{a_2}{b_1b_2}=\frac{c_0c_2}{(c_0-c_1)(c_1-c_2)},\nonumber \\
&\cdots,\ d_n=\frac{a_n}{b_{n-1}{b_n}}=\frac{c_{n-2}c_n}{(c_{n-2}-c_{n-1})(c_{n-1}-c_n)},
\end{align}
for $n\geq 2$.

For the later convenience, let us give general names for two expressions of the continued fraction such that
\begin{align}
&\frac{a_0}{b_0+}\ \frac{a_1}{b_1+}\ \frac{a_2}{b_2+}\ \cdots \equiv
\frac{d_0}{1+}\ \frac{d_1}{1+}\ \frac{d_2}{1+}\ \cdots. 
\end{align}
We call the left {\it standard expression}, and the right {\it normal expression}, while 
the transformation from left to right is unique by
\begin{align}
&\quad d_0=\frac{a_0}{b_0},\ d_1=\frac{a_1}{b_0b_1},\ \cdots,\ d_n=\frac{a_n}{b_{n-1}b_n},\quad (n\geq 1)
\end{align}
but right to left is not unique unless some assumptions are made. \\

\noindent
{\bf (Inversion)}\quad We can express $c$'s in terms of $d$'s inversely such that 
\begin{align}
&c_0=d_0,\ c_1=\frac{d_0d_1}{1+d_1},\ c_2=\frac{d_0d_1d_2}{(1+d_1)(1+d_1+d_2)},\nonumber \\ 
&c_3=\frac{d_0d_1d_2d_3}{(1+d_1+d_2)(1+d_1+d_2+d_3+d_1d_3)},\ \cdots. 
\end{align}
The denominators are expressed by determinants such as
\begin{align}
&D_0\equiv 1,\ D_1\equiv 1+d_1=\left|\begin{array}{cc}
1&d_1\\
-1&1\end{array}\right|,\
D_2\equiv 1+d_1+d_2=\left|\begin{array}{ccc}
1&d_1&0\\
-1&1&d_2\\
0&-1&1\end{array}\right|,
\nonumber \\
&D_3\equiv 1+d_1+d_2+d_3+d_1d_3=\left|\begin{array}{cccc}
1&d_1&0&0\\
-1&1&d_2&0\\
0&-1&1&d_3\\
0&0&-1&1\end{array}\right|,\ \cdots,
\end{align}
whose general expression is obviously given by
\begin{align}
D_n=\left|\begin{array}{ccccc}
1&d_1&0&\cdots&0\\
-1&1&d_2&0&\cdots\\
\cdots&\cdots&\cdots&\cdots&\cdots\\
\cdots&\cdots&-1&1&d_{n}\\
0&\cdots&\cdots&-1&1\end{array}\right|.
\end{align}
By using these results, we have the inversion formula
\begin{align}
c_0=d_0,\quad c_n=\frac{d_0d_1\cdots d_n}{D_{n-1}D_n},\qquad (n\geq 1),
\label{InverseE}
\end{align}
which are new results as far as the author knows. \\

\noindent
{\bf (Examples)}\quad Let us give some examples of Euler type continued fractions. The first example is
\begin{align}
\sum_{n=0}^\infty (-1)^n q^{n(n+1)/2}=\frac{1}{1+}\ \frac{q}{1-q+}\ \frac{q^2}{1-q^2+}\ \frac{q^3}{1-q^3+}\ \cdots,
\label{EulerG_1}
\end{align}
which is derived from another general formula, that is equivalent to Euler formula \eqref{EulerCF},
\begin{align}
a_0+a_0a_1+a_0a_1a_2+\cdots=\frac{a_0}{1-}\ \frac{a_1}{1+a_1-}\ \frac{a_2}{1+a_2-}\ \cdots,
\label{EulerA}
\end{align}
by setting $a_0=1,\ a_1=-q,\ a_2=-q^2,\ \cdots,\ a_n=-q^{n},\ (n\geq 1)$.  

It is instructive to compare \eqref{EulerG_1} with \eqref{Gauss_1} which has the same series in the left hand side,  
implying that continued fraction depends on its computing method, {\it i.e.} non-uniqueness of $C$ part. 

And the second example is
\begin{align}
\sum_{n=0}^\infty (-1)^n q^{n^2}=\frac{1}{1+}\ \frac{q}{1-q+}\ \frac{q^3}{1-q^3+}\ \frac{q^5}{1-q^5+}\ \cdots,
\end{align}
which is derived from \eqref{EulerA} by setting $a_0=1,\ a_1=-q,\ a_2=-q^3,\ \cdots$,\ $a_n=-q^{2n-1},\ (n\geq 1)$. 
Euler gave some other examples in \cite{Euler}, which are omitted here.

\subsection{Gauss-Heine formula for hypergeometric functions}

{\bf (Hypergeometric functions)}\quad 
Hypergeometric functions of Gauss, known earlier also by Euler, are defined by
\begin{align}
{}_2F_1(a, b, c; x)=\sum_{n=0}^\infty \frac{(a)_n(b)_n}{(c)_n n!}\ x^n,
 \end{align}
where $(a)_n$ is the Pochhammer symbol already defined by \eqref{Pochhammer}. 
Its $q$-extension by Heine, probably known also by Gauss, is defined by 
\begin{align}
{}_2\phi_1(\alpha, \beta, \gamma; q, x)=\sum_{n=0}^\infty \frac{(\alpha; q)_n(\beta; q)_n}{(\gamma; q)_n(q; q)_n}\ x^n,
\end{align}
where $(\alpha; q)_n$ is the $q$-Pochhammer symbol also defined by \eqref{qPochhammer}. 
If we set $\alpha=q^a$, we have the limit relation
\begin{align}
\lim_{q\rightarrow 1}\frac{(\alpha; q)_n}{(1-q)^n}=(a)_n,\qquad (\alpha=q^a).
\end{align}
Therefore we have under assumptions $\alpha=q^a,\ \beta=q^b,\ \gamma=q^c$,
\begin{align}
\lim_{q\rightarrow 1} {}_2\phi_1 (\alpha, \beta, \gamma; q, x)= {}_2F_1(a, b, c; x),
\end{align}
which is the reason why ${}_2\phi_1$ is called a $q$-extension of ${}_2F_1$. 

Now let us give {\it contiguous relations} among $q$-hypergeometric functions, 
which can be used for making continued fractions. 
For the sake of simplicity we write in short 
\begin{align}
{}_2\phi_1(\alpha, \beta, \gamma; q, x) \equiv (\alpha,\beta,\gamma; x),
\end{align}
for a while. 
Then contiguous relations in question are written as
\begin{align}
&(1)\quad (\alpha, \beta, \gamma/q; x)-(\alpha, \beta, \gamma; x)
\nonumber \\
&\qquad\qquad =\gamma x\ \frac{(1-\alpha)(1-\beta)}{(q-\gamma)(1-\gamma)}
\cdot(\alpha q,\beta q,\gamma q; x),
\label{ce1}\\
&(2)\quad (\alpha q,\beta,\gamma; x)-(\alpha,\beta,\gamma; x)
\nonumber \\
&\qquad\qquad =\alpha x\ \frac{(1-\beta)}{(1-\gamma)}
\cdot(\alpha q,\beta q,\gamma q; x),\\
&(3)\quad (\alpha q,\beta,\gamma q; x)-(\alpha,\beta,\gamma; x)
\nonumber \\
&\qquad\qquad =x\ \frac{(1-\beta)(\alpha-\gamma)}{(1-\gamma)(1-\gamma q)} 
\cdot(\alpha q,\beta q,\gamma q^2; x),\\
&(4)\quad (\alpha q,\beta/q,\gamma; x)-(\alpha,\beta,\gamma; x)
\nonumber \\
&\qquad\qquad =x\ \frac{(\alpha-\beta/q)}{(1-\gamma)} 
\cdot(\alpha q,\beta,\gamma q; x),
\label{ce4}
\end{align}
which are given, {\it e.g.} in exercise 1.9 of the book by Gasper-Rahman \cite{GasperRahman}, 
although they wrote $\phi(\alpha,\beta,\gamma)$ instead of $(\alpha,\beta,\gamma;x)$. 
For the completeness let us give a derivation of (1)
\begin{align}
&(\alpha, \beta, \gamma/q; x)-(\alpha,\beta,\gamma; x)=
\sum_{n=0}^\infty (\alpha; q)_n(\beta; q)_n \frac{x^n}{(q;q)_n}\left(\frac{1}{(\gamma/q; q)_n}-\frac{1}{(\gamma; q)_n}\right),
\nonumber \\
&\quad=\frac{\gamma/q}{1-\gamma/q}\sum_{n=0}^\infty \frac{(\alpha;q)_n(\beta;q)_n}{(\gamma;q)_n(q;q)_n}x^n(1-q^n)
\nonumber \\
&\quad=\gamma \ \frac{(1-\alpha)(1-\beta)}{(q-\gamma)(1-\gamma)}
\sum_{m=0}^\infty \frac{(\alpha q;q)_m(\beta q;q)_m}{(\gamma q;q)_m(q;q)_m}\ x^{m+1},
\nonumber \\
&\quad=\gamma x\ \frac{(1-\alpha)(1-\beta)}{(q-\gamma)(1-\gamma)}
\cdot(\alpha q,\beta q,\gamma q; x),\nonumber
\end{align}
which is to be derived. Here we set $n=m+1$ and relations $(q;q)_{m+1}=(q;q)_m\cdot(1-q^{m+1})$,  
$(\alpha; q)_{m+1}=(1-\alpha)\cdot(\alpha q;q)_m$, etc. are used. \\

\noindent
{\bf (Continued fractions)}\quad Now let us give an example of making continued fraction, only one from many. 
Using (3) and $(3)'$: an exchange of $\alpha$ and $\beta$ in (3),
\begin{align}
&(3)'\quad (\alpha,\beta q,\gamma q; x)-(\alpha,\beta,\gamma; x)
\nonumber \\
&\qquad\qquad =x\ \frac{(1-\alpha)(\beta-\gamma)}{(1-\gamma)(1-\gamma q)} 
\cdot(\alpha q,\beta q,\gamma q^2; x),
\end{align}
with changing $\alpha, \beta, \gamma \rightarrow \alpha q, \beta, \gamma q$ respectively, 
we have after repetitions
\begin{align}
&\frac{(\alpha q, \beta, \gamma q; x)}{(\alpha,\beta,\gamma; x)}=
\frac{1}{1-}\ \frac{\frac{(1-\beta)(\alpha-\gamma)}{(1-\gamma)(1-\gamma q)}x}
{\frac{(\alpha q,\beta, \gamma q; x)}{(\alpha q,\beta q,\gamma q^2; x)}}
=\frac{1}{1-}\ \frac{\frac{(1-\beta)(\alpha-\gamma)}{(1-\gamma)(1-\gamma q)}x}{1-}\ 
\frac{\frac{(1-\alpha q)(\beta-\gamma q)}{(1-\gamma q)(1-\gamma q^2)}x}
{\frac{(\alpha q,\beta q,\gamma q^2; x)}{(\alpha q^2,\beta q ,\gamma q^3; x)}}
\nonumber \\
&\quad =\frac{1}{1-}\ \frac{\frac{(1-\beta)(\alpha-\gamma)}{(1-\gamma)(1-\gamma q)}x}{1-}\ 
\frac{\frac{(1-\alpha q)(\beta -\gamma q)}{(1-\gamma q)(1-\gamma q^2)}x}{1-}\ 
\frac{\frac{q(1-\beta q)(\alpha-\gamma q)}{(1-\gamma q^2)(1-\gamma q^3)}x}{1-}\ \cdots,
\nonumber \\
&\quad=\frac{e_0}{1-}\ \frac{e_1x}{1-}\ \frac{e_2x}{1-}\ \frac{e_3x}{1-}\ \frac{e_4x}{1-}\ \cdots,
\end{align}
where we set
\begin{align}
&e_0=1,\ e_1=\frac{(1-\beta)(\alpha-\gamma)}{(1-\gamma)(1-\gamma q)},\nonumber \\
&e_2=\frac{(1-\alpha q)(\beta-\gamma q)}{(1-\gamma q)(1-\gamma q^2)},\ 
e_3=q\cdot\frac{(1-\beta q)(\alpha-\gamma q)}{(1-\gamma q^2)(1-\gamma q^3)},\ \cdots,
\label{e_2p1}
\end{align}
whose general expressions are given by 
\begin{align}
\begin{split}
&e_{2m}=q^{m-1}\cdot \frac{(1-\alpha q^m)(\beta-\gamma q^m)}{(1-\gamma q^{2m-1})(1-\gamma q^{2m})},\\ 
&e_{2m+1}=q^m\cdot \frac{(1-\beta q^m)(\alpha-\gamma q^m)}{(1-\gamma q^{2m})(1-\gamma q^{2m+1})},
\end{split}
\end{align}
for $m\geq 1$.\\

\noindent
{\bf (Examples)}\quad Let us give some examples. The first is $\alpha=\beta=0$ case, 
\begin{align}
(\gamma; x)&\equiv {}_0\phi_1(\gamma; q, x)=\sum_{n=0}^\infty \frac{x^n}{(\gamma;q)_n(q; q)_n}
\\
&=1+\frac{x}{(1-\gamma)(1-q)}+\frac{x^2}{(1-\gamma)(1-\gamma q)(1-q)(1-q^2)}+\cdots.
\nonumber
\end{align}
The continued fraction considered here is not 
for $(\gamma; x)$ itself, but for the quotient
\begin{align}
\frac{(\gamma q; x)}{(\gamma; x)}=\frac{e_0}{1-}\ \frac{e_1 x}{1-}\ \frac{e_2 x}{1-}\ 
\frac{e_3 x}{1-}\ \cdots,
\end{align}
where the coefficients are given by \eqref{e_2p1} with $\alpha=\beta=0$, that is,
\begin{align}
\begin{split}
&e_0=1,\ e_1=-\frac{\gamma}{(1-\gamma)(1-\gamma q)},
\\ 
&e_2=-\frac{\gamma q}{(1-\gamma q)(1-\gamma q^2)},\ e_3=-\frac{\gamma q^2}{(1-\gamma q^2)(1-\gamma q^3)},\ \cdots,
\end{split}
\end{align}
general rule of which will be obvious.

The second example is $\alpha=0$ case, that is
\begin{align}
(\beta, \gamma; x)&\equiv {}_1\phi_1(\beta,\gamma; q, x)=\sum_{n=0}^\infty \frac{(\beta;q)_n}{(\gamma; q)_n}\ 
\frac{x^n}{(q;q)_n}\\
&=1+\frac{(1-\beta)x}{(1-\gamma)(1-q)}+\frac{(1-\beta)(1-\beta q)x^2}{(1-\gamma)(1-\gamma q)(1-q)(1-q^2)}+\cdots.
\nonumber
\end{align}
It should be noted that there are two types of quotient 
\begin{align}
\frac{(\beta,\gamma q;x)}{(\beta,\gamma;x)}\quad \text{and}\quad  
\frac{(\beta q,\gamma q;x)}{(\beta,\gamma;x)},
\nonumber
\end{align}
whose continued fractions are a little bit different each other. 
The continued fraction for the first case is simply given by setting $\alpha=0$ in \eqref{e_2p1}.

The result for the second case is obtained as follows. Due to the symmetry between $\alpha$ and $\beta$, 
exchange the names of $\alpha$ and $\beta$ first then after that set $\alpha=0$ in \eqref{e_2p1}. 
Such procedure of exchanging $\alpha$ and $\beta$ corresponds to the order of applying contiguous relations 
$(3)'$ after (3) or (3) after $(3)'$.
To make sure, let us write down the result explicitly,
\begin{align}
\frac{(\beta q, \gamma q; x)}{(\beta, \gamma; x)}=\frac{e_0}{1-}\ \frac{e_1 x}{1-}\ \frac{e_2 x}{1-}\ 
\frac{e_3 x}{1-}\ \cdots,
\end{align}
where the coefficients are given by
\begin{align}
\begin{split}
&e_0=1,\ e_1=-\frac{\gamma-\beta}{(1-\gamma)(1-\gamma q)},
\\
&e_2=-\frac{\gamma q(1-\beta q)}{(1-\gamma q)(1-\gamma q^2)},\ e_3=-\frac{q(\gamma q-\beta)}{(1-\gamma q^2)(1-\gamma q^3)},\ 
\cdots,
\end{split}
\end{align}
whose general expressions are written by
\begin{align}
\begin{split}
e_{2m}=-\frac{\gamma q^{2m-1}(1-\beta q^m)}{(1-\gamma q^{2m-1})(1-\gamma q^{2m})},
\\
e_{2m+1}=-\frac{q^m(\gamma q^m-\beta)}{(1-\gamma q^{2m})(1-\gamma q^{2m+1})},
\end{split}
\end{align}
for $m\geq 1$. 

\subsection{General formulas by Muir and Rogers}

\noindent
{\bf (Muir formula)}\quad 
Muir \cite{Muir} derived in 1875 his general formula transforming a quotient of two infinite power series 
into a continued fraction
\begin{align}
\frac{c_0+c_1x+c_2x^2+\cdots}{b_0+b_1x+b_2x^2+\cdots}=
\frac{e_0}{1-}\ \frac{e_1x}{1-}\ \frac{e_2x}{1-}\ \cdots,
\label{MuirSC}
\end{align}
where $x$ plays the role of pivot variable. 
Following Muir, let us introduce determinants defined by
\begin{align}
&\theta_0=b_0,\ \theta_1=c_0,\  
\theta_2=\left|\begin{array}{cc}
b_0&b_1\\
c_0&c_1\end{array}\right|,\ 
\theta_3=\left|\begin{array}{ccc}
b_0&b_1&b_2\\
0&c_0&c_1\\
c_0&c_1&c_2\end{array}\right|,
\nonumber \\
&\theta_4=\left|\begin{array}{cccc}
b_0&b_1&b_2&b_3\\
0&b_0&b_1&b_2\\
0&c_0&c_1&c_2\\
c_0&c_1&c_2&c_3\end{array}\right|,\ 
\theta_5=\left|\begin{array}{ccccc}
b_0&b_1&b_2&b_3&b_4\\
0&b_0&b_1&b_2&b_3\\
0&0&c_0&c_1&c_2\\
0&c_0&c_1&c_2&c_3\\
c_0&c_1&c_2&c_3&c_4\end{array}\right|,\ \cdots
\label{tMuir}
\end{align}
whose construction rule will be almost obvious; that is, 
$\theta_{2m}$ consists of the same number of rows given by $b$'s and $c$'s, and $\theta_{2m+1}$ has 
$m$ rows of $b$'s and $m+1$ rows of $c$'s, having a special configurational structure.

Then the coefficients $e_0,\ e_1,\ e_2,\ \cdots$, in \eqref{MuirSC} are given by
\begin{align}
&e_0=\frac{\theta_1}{\theta_0},\ 
e_1=\frac{\theta_2}{\theta_1\theta_0},\ 
e_2=\frac{\theta_3}{\theta_2\theta_1},\ 
e_3=\frac{\theta_4\theta_1}{\theta_3\theta_2},\ \cdots,
\nonumber \\
&e_n=\frac{\theta_{n+1}\theta_{n-2}}{\theta_{n}\theta_{n-1}},\quad (n\geq 3),
\label{MuirF}
\end{align}
which we call Muir formula.

It should be noted that Muir formula contains previous Gauss-Heine formula as a special case whose coefficients $b$'s and $c$'s 
are connected in a special fashion. 
General Muir formula can be applied to cases having no relation among $b$'s and $c$'s, however the rule of making coefficients 
$e$'s may become complicated for such cases. \\

\noindent
{\bf (Example)}\quad Let us apply Muir formula \eqref{MuirF} to Rogers-Ramanujan quotient defined by 
\begin{align}
\frac{F(qx)}{F(x)}=\frac{e_0}{1-}\ \frac{e_1x}{1-}\ \frac{e_2x}{1-}\ \cdots,
\end{align}
where we introduced
\begin{align}
F(x)=1+\sum_{n=1}^\infty \frac{q^{n^2}}{(q;q)_n}x^n,\quad
F(qx)=1+\sum_{n=1}^\infty \frac{q^{n(n+1)}}{(q;q)_n}x^n.
\end{align}
The coefficients are given by $b_0=1,\ c_0=1$ and
\begin{align}
\begin{split}
&b_1=\frac{q}{1-q},\ b_2=\frac{q^4}{(1-q)(1-q^2)},\ b_3=\frac{q^9}{(1-q)(1-q^2)(1-q^3)},\ \cdots, \\
&c_1=\frac{q^2}{1-q},\ c_2=\frac{q^6}{(1-q)(1-q^2)},\ c_3=\frac{q^{12}}{(1-q)(1-q^2)(1-q^3)},\ \cdots,
\end{split}
\nonumber
\end{align}
then from \eqref{tMuir} we have $\theta$'s such that
\begin{align}
\theta_0=1,\ \theta_1=1,\ \theta_2=-q,\ \theta_3=q^3,\ \theta_4=q^7,\ \theta_5=q^{13},\ \cdots,
\nonumber
\end{align}
from which we obtain by using \eqref{MuirF},
\begin{align}
e_0=1,\ e_1=-q,\ e_2=-q^2,\ e_3=-q^3,\ e_4=-q^4,\ \cdots.
\end{align}
This is the famous Rogers-Ramanujan identity
\begin{align}
\frac{F(qx)}{F(x)}=\frac{1}{1+}\ \frac{qx}{1+}\ \frac{q^2x}{1+}\ \frac{q^3x}{1+}\ \cdots,
\label{RRM}
\end{align}
where usually $x=1$ is assumed. Another derivation will be given in \S 2.4 later.\\

\noindent
{\bf (Muir-Rogers formula)}\quad 
For the case of single series (not quotient), that is, $b_0=1,\ b_1=b_2=\cdots=0$, 
was discussed also by Rogers \cite{Rogers} 
who introduced determinants
\begin{align}
&\alpha_0=c_0,\ \alpha_1=c_1,\ \alpha_2=\left|\begin{array}{cc}
c_0&c_1\\
c_1&c_2\end{array}\right|,\ 
\alpha_3=\left|\begin{array}{cc}
c_1&c_2\\
c_2&c_3\end{array}\right|,
\nonumber \\
&\alpha_4=\left|\begin{array}{ccc}
c_0&c_1&c_2\\
c_1&c_2&c_3\\
c_2&c_3&c_4\end{array}\right|,\ 
\alpha_5=\left|\begin{array}{ccc}
c_1&c_2&c_3\\
c_2&c_3&c_4\\
c_3&c_4&c_5\end{array}\right|,\ \cdots,
\label{aRogers}
\end{align}
which are known as Hankel type determinants.

Determinants of the Muir formula \eqref{tMuir} become \eqref{aRogers}, because by setting $b_0=1,\ b_1=b_2=\cdots=0$, we obtain    
$\theta_0=1,\ \theta_1=c_0=\alpha_0,\ \theta_2=c_1=\alpha_1,\ \theta_3=\alpha_2,\ \theta_4=\alpha_3$, and so on.  
Since we have $\theta_{n+1}=\alpha_n,\ (n\geq 0)$ in general, Muir formula \eqref{MuirF} becomes
\begin{align}
&e_0=\alpha_0,\ e_1=\frac{\alpha_1}{\alpha_0},\ e_2=\frac{\alpha_2}{\alpha_1\alpha_0},\ 
e_3=\frac{\alpha_3\alpha_0}{\alpha_2\alpha_1},\ \cdots
\nonumber \\
&e_n=\frac{\alpha_n\alpha_{n-3}}{\alpha_{n-1}\alpha_{n-2}},\quad (n\geq 3),
\label{MuirRogersF}
\end{align}
which we call Muir-Rogers formula (or simply Rogers formula). \\

\noindent
{\bf (Inversion)}\quad The above Muir and Muir-Rogers formulas are mappings from quotient of series or single series ($S$) 
to continued fraction ($C$). 
Now let us consider about their inverse mapping. 
Since the former (Muir formula) is two to one correspondence, and the inversion is not unique, we restrict ourselves 
only to Muir-Rogers formula \eqref{MuirRogersF}, that is, we want to find a mapping from $e$'s to $c$'s 
with intermediate $\alpha$'s.  

Such procedure is described as follows. 
The first step is from $e$'s to $\alpha$'s. Since we have $\alpha_0=e_0=c_0,\ \alpha_1=e_0e_1=c_1$, and
\begin{align}
e_0e_1\cdots e_{2m}=\frac{\alpha_{2m}}{\alpha_{2m-2}},\quad e_0e_1\cdots e_{2m+1}=\frac{\alpha_{2m+1}}{\alpha_{2m-1}},
\end{align}
from which we obtain
\begin{align}
\begin{split}
&\alpha_{2m}=(e_{2m}e_{2m-1})\cdot(e_{2m-2}e_{2m-3})^2\cdots (e_2e_1)^m \cdot (e_0)^{m+1},\\
&\alpha_{2m+1}=(e_{2m+1}e_{2m})\cdot(e_{2m-1}e_{2m-2})^2\cdots (e_3e_2)^m\cdot(e_1e_0)^{m+1}.
\end{split}
\end{align}
The second step is from determinants $\alpha$'s to matrix elements $c$'s, that is, to find 
$c_0=\alpha_0,\ c_1=\alpha_1$ and find $c_2,\ c_3$ from
\begin{align}
&\left|\begin{array}{cc}
c_0&c_1\\
c_1&c_2\end{array}\right|=c_0c_2-c_1^2=\alpha_2,\ 
\left|\begin{array}{cc}
c_1&c_2\\
c_2&c_3\end{array}\right|=c_1c_3-c_2^2=\alpha_3,
\label{InverseMR}
\end{align}
and so on sequentially, which are tedious but straight-forward tasks. 
We will give later in \S 3.1 one example of such inversion calculation. \\

\noindent
{\bf (Examples)}\quad Let us give some examples of Muir-Rogers formula. 
As was mentioned before the Cauchy's identity \eqref{Cauchy} is such an example, 
\begin{align}
{}_1\phi_0(\alpha; q, x)=\sum_{n=0}^\infty \frac{(\alpha; q)_n}{(q; q)_n}\ x^n=\frac{e_0}{1-}\ \frac{e_1x}{1-}\ \frac{e_2x}{1-}\ \cdots.
\end{align}
Since $c_n=(\alpha; q)_n/(q; q)_n$ we have $\alpha_0=e_0=c_0=1$ and 
\begin{align}
\alpha_1=c_1=\frac{1-\alpha}{1-q},\quad 
\alpha_2=\left|\begin{array}{cc}
c_0&c_1\\
c_1&c_2\end{array}\right|=\frac{(1-\alpha)(\alpha-q)}{(1-q)(1-q^2)},\ \cdots,
\end{align}
thus we obtain
\begin{align}
&e_0=\alpha_0=1,\ e_1=\frac{\alpha_1}{\alpha_0}=\frac{1-\alpha}{1-q},
\nonumber \\
&e_2=\frac{\alpha_2}{\alpha_1\alpha_0}=\frac{\alpha-q}{1-q^2},\ 
e_3=\frac{\alpha_3\alpha_0}{\alpha_2\alpha_1}=\frac{q(1-\alpha q)}{(1+q)(1-q^3)},\ \cdots,
\nonumber \\
&e_{2m}=\frac{q^{m-1}(\alpha-q^m)}{(1+q^m)(1-q^{2m-1})},\ 
e_{2m+1}=\frac{q^m(1-\alpha q^m)}{(1+q^m)(1-q^{2m+1})},
\end{align}
for $m\geq 1$, which are to be derived. 

Similarly let us give the next example
\begin{align}
&{}_1\rho_0(\alpha; x)=\sum_{n=0}^\infty (\alpha; q)_n\ x^n=\frac{e_0}{1-}\ \frac{e_1x}{1-}\ \frac{e_2x}{1-}\ \cdots,
\end{align}
whose coefficients are given by
\begin{align}
&e_0=1,\ e_1=1-\alpha,\ e_2=\alpha(1-q),\ e_3=q(1-\alpha q),\ \cdots.
\nonumber \\
&e_{2m}=\alpha q^{m-1}(1-q^m),\ e_{2m+1}=q^m(1-\alpha q^m),\quad (m\geq 1).
\end{align}
This series is $q$-extension of the Wallis series \eqref{WF} (with $x\rightarrow -x$), and 
the case $\alpha=q$ and taking limit $q\rightarrow 1$ with $x=-1/(1-q)$ is the original 
Wallis series \eqref{WallisSum}.

\subsection{Ramanujan formula for his continued fractions}

\noindent
{\bf (Ramanujan formula)}\quad 
According to Andrews \cite{Andrews2}, Ramanujan considered a limit of Heine's $q$-hypergeometric function 
\begin{align}
{}_2\phi_1(\alpha, \beta, \gamma; q, z)=\sum_{n=0}^\infty \frac{(\alpha; q)_n(\beta; q)_n}{(\gamma; q)_n(q; q)_n}\ z^n,
\end{align}
such that $\alpha\rightarrow\infty,\ z\rightarrow 0$ with a constraint $\alpha z=-qx$, which becomes
\begin{align}
{}_1\Phi_1(\beta, \gamma; q, x)=\sum_{n=0}^\infty q^{n(n+1)/2}\frac{(\beta; q)_n}{(\gamma; q)_n}\ \frac{x^n}{(q; q)_n}.
\end{align}
Starting from this, he introduced sequence of functions ($k=0, 1, 2, \cdots$)
\begin{align}
&F_k(\beta, \gamma; q, x)=\sum_{n=0}^\infty q^{n(n+2k+1)/2}\frac{(\beta q^k; q)_n}{(\gamma; q)_{n+k}}\ \frac{x^n}{(q; q)_n},
\\
&H_k(\beta, \gamma; q, x)=\sum_{n=0}^\infty q^{n(n+2k+3)/2}\frac{(\beta q^k; q)_n}{(\gamma; q)_{n+k}}\ \frac{x^n}{(q; q)_n}.
\end{align}
It should be noted that they satisfy relations
\begin{align}
&F_0(\beta, \gamma; q; x)={}_1\Phi_1(\beta, \gamma; q, x),\\ 
&H_k(\beta, \gamma; q, x)=F_k(\beta, \gamma; q, qx).
\end{align}
Here we have changed notations from those of Andrews (or Ramanujan): $\beta=-\lambda/a,\ \gamma=-bq$, and $x=a$ for 
convenience, we keep our notations for a while.  

To derive Ramanujan formula for continued fractions, we need the contiguous relations among $F$'s and $H$'s: 
\begin{align}
(1)\quad F_k(\beta, \gamma; q, x)&-H_k(\beta, \gamma; q, x)
\nonumber \\
&=\left(q^{k+1}-\beta q^{2k+1}\right) x\cdot F_{k+1}(\beta, \gamma; q, x),
\label{R1}
\\
(2)\quad H_k(\beta, \gamma; q, x)&-F_{k+1}(\beta, \gamma; q, x)
\nonumber \\
&=-\left(\gamma q^k+\beta q^{2k+2}x\right)\cdot H_{k+1}(\beta, \gamma; q, x),
\label{R2} 
\end{align}
where one should notice the positions of pivot variable $x$. 
Except $\gamma=0$ case the above one does not belong to Gauss-Heine type contiguous relations. 
According to Ramanujan's notation, the coefficients in the right hand side become
\begin{align}
\begin{split}
&\left(q^{k+1}-\beta q^{2k+1}\right) x=aq^{k+1}+\lambda q^{2k+1},\\ 
&-\left(\gamma q^k+\beta q^{2k+2}x\right)=bq^{k+1}+\lambda q^{2k+2}, 
\end{split}
\end{align}
which agree with those of Andrews \cite{Andrews2}.

For the sake of completeness, let us show \eqref{R1} and \eqref{R2} by direct substitutions, contrary to method of 
Andrews who used Heine's transformation relation. 
The first one is derived as follows. 
\begin{align}
&F_k-H_k=\sum_{n=0}^\infty q^{n(n+2k+1)/2} \frac{(\beta q^k; q)_n}{(\gamma;q)_{n+k}}
\frac{x^n}{(q; q)_n} \left(1-q^n\right)
\nonumber \\
&\ =\sum_{m=0}^\infty q^{m(m+2k+3)/2}q^{k+1}\frac{(\beta q^{k}; q)_{m+1}}{(\gamma; q)_{m+k+1}}
\frac{x^{m+1}}{(q;q)_{m+1}}
(1-q^{m+1}),
\label{part1} 
\end{align}
where we set $n=m+1$, since $n=0$ term vanishes, and use equality $n(n+2k+1)/2=m(m+2k+3)/2+(k+1)$. Then by using 
$(q;q)_{m+1}=(q;q)_m\cdot(1-q^{m+1})$ and $(\beta q^k; q)_{m+1}=(1-\beta q^k)\cdot(\beta q^{k+1}; q)_m$, 
the right hand side (RHS1) of \eqref{part1} becomes
\begin{align}
\text{RHS1}&=q^{k+1}\left(1-\beta q^k\right) x  
\sum_{m=0}^\infty q^{m(m+2(k+1)+1)/2}\ \frac{(\beta q^{k+1}; q)_m}{(\gamma; q)_{m+k+1}}\ \frac{x^m}{(q;q)_m}
\nonumber \\
&=\left(q^{k+1}-\beta q^{2k+1} \right) x\cdot F_{k+1}(\beta, \gamma; q, x),
\end{align}
which is to be derived. The second one is derived as follows. 
\begin{align}
&H_k-F_{k+1}
\nonumber \\
&\quad=\sum_{n=0}^\infty q^{n(n+2k+3)/2}\frac{x^n}{(q;q)_n}\left(
\frac{(\beta q^k; q)_n}{(\gamma; q)_{n+k}}-\frac{(\beta q^{k+1}; q)_n}{(\gamma; q)_{n+k+1}}\right),
\label{part2}
\end{align}
where the parentheses part is rewritten by
\begin{align} 
\left(\cdots\right)&=\frac{(\beta q^{k+1}; q)_n}{(\gamma; q)_{n+k+1}}\left[
(1-\gamma q^{n+k})\left(\frac{1-\beta q^k}{1-\beta q^{k+n}}\right)-1\right]
\nonumber \\
&=\frac{(\beta q^{k+1}; q)_n}{(\gamma; q)_{n+k+1}}
\left[ -\gamma q^{k+n}-\beta q^k\ \frac{1-\gamma q^{n+k}}{1-\beta q^{k+n}}\ (1-q^n) \right].
\nonumber
\end{align}
Therefore we obtain the right hand side (RHS2) of \eqref{part2}, by changing $n=m+1$ as before in the second term, 
\begin{align}
\text{RHS2}&=-\gamma q^k \sum_{n=0}^\infty q^{n(n+2k+5)/2} \frac{x^n}{(q; q)_n}\cdot\frac{(\beta q^{k+1}; q)_n}{(\gamma; q)_{n+k+1}}
\nonumber \\
&\quad -\beta q^{2k+2} \sum_{m=0}^\infty q^{m(m+2k+5)/2} \frac{x^{m+1}}{(q; q)_m}\cdot\frac{(\beta q^{k+1}; q)_m}{(\gamma; q)_{m+k+1}}
\nonumber \\
&=-\left(\gamma q^k+\beta q^{2k+2} x \right)\cdot H_{k+1}(\beta, \gamma; q, x),
\end{align}
which is to be derived.

By use of contiguous relations \eqref{R1} and \eqref{R2}, we obtain the continued fraction
\begin{align}
&\frac{H_0}{F_0}=\frac{1}{1+}\ \frac{(q-\beta q)x}{1-}\ \frac{(\gamma+\beta q^2 x)}{F_1/H_1}
\nonumber \\
&=\frac{1}{1+}\ \frac{(q-\beta q)x}{1-}\ \frac{(\gamma+\beta q^2 x)}{1+}\ \frac{(q^2-\beta q^3)x}{1-}\ 
\frac{(\gamma q +\beta q^4 x)}{1+}\ \cdots,
\label{MR_F}
\\
&=\frac{1}{1+}\ \frac{aq+\lambda q}{1+}\ \frac{bq+\lambda q^2}{1+}\ \frac{aq^2+\lambda q^3}{1+}\ 
\frac{bq^2+\lambda q^4}{1+}\ \cdots,
\label{RamanujanF}
\end{align}
the latter of which is written in Ramanujan's notation ($\beta=-\lambda/a,\ \gamma=-bq,\ x=a$).
This is the Ramanujan's continued fraction to be derived. 

Ramanujan uses \eqref{RamanujanF} {\it reversely}, by choosing values of $a, \lambda$, and  $b$ variously to find 
seiries $F_0,\ H_0$ which he renamed as
\begin{align}
\begin{split}
&F_0(\beta, \gamma; q, x)=G(a, \lambda; b, q),\\
&H_0(\beta, \gamma; q, x)=G(aq, \lambda q; b, q),
\end{split}
\end{align}
which are given explicitly by
\begin{align}
\begin{split}
&G(a, \lambda; b, q)=1+\sum_{n=1}^\infty \frac{q^{n(n+1)/2}}{(q;q)_n}\cdot
\frac{(a+\lambda)\cdots(a+\lambda q^{n-1})}{(1+bq)\cdots(1+bq^n)},\\
&G(aq, \lambda q; b, q)=1+\sum_{n=1}^\infty \frac{q^{n(n+3)/2}}{(q;q)_n}\cdot
\frac{(a+\lambda)\cdots(a+\lambda q^{n-1})}{(1+bq)\cdots(1+bq^n)},
\end{split}
\end{align}
where $\beta=-\lambda/a,\ \gamma=-bq,\ x=a$ are used as before, \\

\noindent
{\bf (Example)}\quad Let us give here only one example. It is the Rogers-Ramanujan identity, 
whose continued fraction is given simply by setting $a=0,\ \lambda=1,\ b=0$, which gives
\begin{align}
\begin{split}
&F_0=G(0, 1; 0, q)=\sum_{n=0}^\infty \frac{q^{n^2}}{(q; q)_n},\\
&H_0=G(0, q; 0, q)=\sum_{n=0}^\infty \frac{q^{n(n+1)}}{(q; q)_n},
\end{split}
\end{align}
which is the known result, while if $\lambda$ is kept as an independent parameter $x$ it can be 
regarded as the pivot variable, as was shown in \eqref{RRM} before. This is possible because $b=0$, {\it i.e.} $\gamma=0$.

Other examples will be given in \S 3.1 later, which are taken from examples given by Andrews \cite{Andrews2}. 

\subsection{Gauss-Heine type formula for extended hypergeometric functions}

\noindent
{\bf (Muir type continued fractions)}\quad 
In the previous subsection we considered a quotient of two series
\begin{align}
\begin{split}
&F_0(\beta, \gamma; q, x)=\sum_{n=0}^\infty q^{n(n+1)/2}\ \frac{(\beta; q)_n}{(\gamma; q)_n}\ \frac{x^n}{(q;q)_n},\\
&H_0(\beta, \gamma; q, x)=\sum_{n=0}^\infty q^{n(n+3)/2}\ \frac{(\beta; q)_n}{(\gamma; q)_n}\ \frac{x^n}{(q;q)_n},
\end{split}
\end{align}
which are both extended $q$-hypergeometric functions. Let us try to apply Muir formula \eqref{MuirF} to $H_0/F_0$ by setting 
coefficients 
\begin{align}
&b_0=1,\ b_1=\frac{q}{1-q}\frac{1-\beta}{1-\gamma},\ 
b_2=\frac{q^3}{(1-q)(1-q^2)}\frac{(1-\beta)(1-\beta q)}{(1-\gamma)(1-\gamma q)},\ \cdots
\nonumber \\
&c_0=1,\ c_1=\frac{q^2}{1-q}\frac{1-\beta}{1-\gamma},\ 
c_2=\frac{q^5}{(1-q)(1-q^2)}\frac{(1-\beta)(1-\beta q)}{(1-\gamma)(1-\gamma q)},\ \cdots
\nonumber
\end{align}
and taking $x$ as the pivot variable. Then we obtain a continued fraction
\begin{align}
&\frac{H_0}{F_0}=\frac{e_0}{1-}\ \frac{e_1x}{1-}\ \frac{e_2x}{1-}\ \frac{e_3x}{1-}\ \cdots,
\label{RMC}
\end{align}
with coefficients given by
\begin{align}
&e_0=1,\ e_1=-\frac{q(1-\beta)}{1-\gamma},\ e_2=-\frac{q^2(\gamma-\beta)}{(1-\gamma)(1-\gamma q)},
\nonumber \\ 
&e_3=-\frac{q^2(1-\beta q)}{(1-\gamma q)(1-\gamma q^2)},\ 
e_4=-\frac{q^4(\gamma q-\beta)}{(1-\gamma q^2)(1-\gamma q^3)},\ \cdots,
\label{Rbc}
\end{align}
where the general coefficients are written by
\begin{align}
\begin{split}
&e_{2m}=-\frac{q^{2m}(\gamma q^{m-1}-\beta)}{(1-\gamma q^{2m-2})(1-\gamma q^{2m-1})},\\
&e_{2m+1}=-\frac{q^{m+1}(1-\beta q^m)}{(1-\gamma q^{2m-1})(1-\gamma q^{2m})},
\end{split}
\end{align}
for $m\geq 1$. 

Note the different position of pivot variable $x$ in \eqref{RMC} compared with \eqref{MR_F}. 
Thus we have another continued fraction for the same $H_0/F_0$, which implies again non-uniqueness of $C$. \\

\noindent
{\bf (Extended contiguous relations)}\quad
Such regular structure suggests us that for extended $q$-hypergeometric functions we can 
expect Gauss-Heine type contiguous relations as in \S2.2. 
Let us consider such extended $q$-hypergeometric 
functions
\begin{align}
{}_2\Phi_1(\alpha, \beta, \gamma; q, x)=\sum_{n=0}^\infty 
q^{n(n+1)/2}\frac{(\alpha; q)_n(\beta; q)_n}{(\gamma; q)_n(q;q)_n}\ x^n,
\end{align}
which becomes previous $F_0(\beta, \gamma; q, x)={}_1\Phi_1(\beta,\gamma;q,x)$ if we set $\alpha=0$. 
Further this is another $q$-extension together with Heine's ${}_2\phi_1$ that becomes Gauss's hypergeometric function 
${}_2F_1$ when $q\rightarrow 1$ limit. 

For the simplicity let us write $\Phi$'s in short by
\begin{align}
{}_2\Phi_1^{(k)}(\alpha,\beta,\gamma;q,x)&=
\sum_{n=0}^\infty q^{n(n+2k+1)/2}\frac{(\alpha;q)_n(\beta; q)_n}{(\gamma;q)_n(q;q)_n}\ x^n
\nonumber \\
&=[\alpha,\beta,\gamma;x]_k,
\end{align}
for a while. Here we set ${}_2\Phi_1^{(k)}(x)={}_2\Phi_1(xq^k),\ (k=0,1,2,\cdots)$. 
Then the contiguous relations analogous to \eqref{ce1} $\sim$ \eqref{ce4} are given by
\begin{align}
&(1)\quad [\alpha,\beta,\gamma/q;x]_k-[\alpha,\beta,\gamma;x]_k
\nonumber \\
&\qquad\qquad=\gamma q^{k+1} x\frac{(1-\alpha)(1-\beta)}{(q-\gamma)(1-\gamma)}\cdot[\alpha q,\beta q,\gamma q;x]_{k+1},
\label{CE1}\\
&(2)\quad [\alpha q,\beta, \gamma; x]_k-[\alpha,\beta,\gamma;x]_k
\nonumber \\
&\qquad\qquad=\alpha q^{k+1}x\ \frac{1-\beta}{1-\gamma}\cdot[\alpha q,\beta q,\gamma q; x]_{k+1},
\\
&(3)\quad [\alpha q,\beta,\gamma q;x]_k-[\alpha,\beta,\gamma;x]_k
\nonumber \\
&\qquad\qquad=q^{k+1}x\ \frac{(1-\beta)(\alpha-\gamma)}{(1-\gamma)(1-\gamma q)}\cdot[\alpha q,\beta q,\gamma q^2; x]_{k+1},
\\
&(4)\quad [\alpha q,\beta/q,\gamma;x]_k-[\alpha,\beta,\gamma;x]_k
\nonumber \\
&\qquad\qquad=q^{k+1}x\ \frac{\alpha-\beta/q}{1-\gamma}\cdot[\alpha q,\beta,\gamma q; x]_{k+1}.
\label{CE4}
\end{align}
Compare these with those of \S 2.2, which are different not only of the existence of $k$, but also a little bit different, 
because the above contiguous relations have additional $q$'s in various places.  
Due to the exchange symmetry between $\alpha$ and $\beta$, we have $(2)'$, $(3)'$ and $(4)'$ as before. 

For the completeness let us derive the relation (1) for an example.
\begin{align}
&[\alpha,\beta,\gamma/q;x]_k-[\alpha,\beta,\gamma;x]_k
\nonumber \\
& =\sum_{n=0}^\infty q^{n(n+2k+1)/2} \frac{(\alpha;q)_n(\beta;q)_n}{(q;q)_n}x^n
\left(\frac{1}{(\gamma/q;q)_n}-\frac{1}{(\gamma;q)_n}\right),
\label{RR1_RHS}
\end{align}
where the parentheses part is rewritten as
\begin{align}
\left( \cdots \right)&=\frac{1}{(1-\gamma)\cdots(1-\gamma q^{n-2})}\left(\frac{1}{1-\gamma/q}-\frac{1}{1-\gamma q^{n-1}}\right)
\nonumber \\
&=\frac{\gamma}{q-\gamma}\frac{1-q^n}{(\gamma;q)_n}, 
\nonumber
\end{align}
therefore we have the right hand side (RHS3) of \eqref{RR1_RHS}
\begin{align}
&\text{RHS3}=\frac{\gamma}{q-\gamma}\sum_{n=0}^\infty q^{n(n+2k+1)/2}\frac{(\alpha;q)_n(\beta;q)_n}{(\gamma;q)_n (q;q)_n}\ 
x^n (1-q^n)
\nonumber \\
&=\frac{\gamma}{q-\gamma}\sum_{m=0}^\infty q^{m(m+2k+3)/2} q^{k+1} 
\frac{(\alpha;q)_{m+1}(\beta;q)_{m+1}}{(\gamma;q)_{m+1}(q;q)_{m+1}} x^{m+1}(1-q^{m+1})
\nonumber \\
&=\frac{\gamma}{q-\gamma}\frac{(1-\alpha)(1-\beta)}{(1-\gamma)}q^{k+1}x \cdot [\alpha q,\beta q,\gamma q; x]_{k+1},
\end{align}
which is to be shown. 

These four contiguous relations however are not enough to make continued fractions, because of the parameter $k$, 
which will be discussed next. \\

\noindent
{\bf (Absence of extended Ramanujan type continued fractions)}\quad 
In addition to the above four relations, we have
\begin{align}
&(0)\quad 
[\alpha,\beta,\gamma;x]_k-[\alpha,\beta,\gamma; x]_{k+1}
\nonumber \\
&\qquad\qquad =q^{k+1}x\ \frac{(1-\alpha)(1-\beta)}{(1-\gamma)}\cdot [\alpha q,\beta q,\gamma q; x]_{k+1},
\label{CE0}
\end{align}
which corresponds to $F_0-H_0=qx(1-\beta)F_{1}$ when $\alpha=0, k=0$. This relation is directly derived such as
\begin{align}
&[\alpha,\beta,\gamma;x]_k-[\alpha,\beta,\gamma; x]_{k+1}
\nonumber \\
&=\sum_{n=0}^\infty q^{n(n+2k+1)/2}\frac{(\alpha;q)_n(\beta;q)_n}{(\gamma;q)_n(q;q)_n}x^n(1-q^n)
\nonumber \\
&=q^{k+1}x\ \frac{(1-\alpha)(1-\beta)}{(1-\gamma)}\cdot [\alpha q,\beta q,\gamma q; x]_{k+1},
\nonumber
\end{align}
which is to be shown.

Now let us set for the present three parameters $(\alpha,\beta,\gamma)$ case also 
\begin{align}
\begin{split}
&\tilde{F}_k(\alpha, \beta, \gamma; x)=\sum_{n=0}^\infty q^{n(n+2k+1)/2}\
\frac{(\alpha q^k; q)_n(\beta q^k;q)_n)}{(\gamma;q)_{n+k}} \frac{x^n}{(q;q)_n},
\\
&\tilde{H}_k(\alpha, \beta, \gamma; x)=\sum_{n=0}^\infty q^{n(n+2k+3)/2}\
\frac{(\alpha q^k; q)_n(\beta q^k;q)_n)}{(\gamma;q)_{n+k}} \frac{x^n}{(q;q)_n},
\end{split}
\end{align}
we can derive similarly the first relation
\begin{align}
&(1)\quad \tilde{F}_k(\alpha,\beta, \gamma; x)-\tilde{H}_k(\alpha,\beta,\gamma; x)
\nonumber \\
&\qquad =q^{k+1}(1-\alpha q^k)(1-\beta q^k)x\cdot\tilde{F}_{k+1}(\alpha, \beta, \gamma; x),
\label{RR1}
\end{align}
which is an extension of \eqref{R1}. 

The second relation expected shall be $\tilde{H}_k-\tilde{F}_{k+1}\propto \tilde{H}_{k+1}$, 
whose left hand side is written by
\begin{align}
&(2)\quad \tilde{H}_k-\tilde{F}_{k+1}=\sum_{n=0}^\infty q^{n(n+2k+3)/2}\frac{x^n}{(q;q)_n}
\nonumber \\
&\qquad\times
\left(\frac{(\alpha q^k;q)_n(\beta q^k;q)_n}{(\gamma;q)_{n+k}}-
\frac{(\alpha q^{k+1};q)_n(\beta q^{k+1};q)_n}{(\gamma;q)_{n+k+1}}\right),
\label{RR2}
\end{align}
where the parentheses part is rewritten by
\begin{align}
\left(\cdots\right)&=\frac{(\alpha q^{k+1};q)_n(\beta q^{k+1};q)_n}{(\gamma;q)_{n+k+1}}
\nonumber \\
&\quad\times
\left[(1-\gamma q^{n+k})\frac{(1-\alpha q^k)}{(1-\alpha q^{k+n})}\frac{(1-\beta q^k)}{(1-\beta q^{k+n})}-1\right],
\nonumber
\end{align}
while 
\begin{align}
\left[\cdots\right]&=-\gamma q^{n+k}-
(1-q^n)(1-\gamma q^{n+k})\left(\frac{\alpha q^k}{1-\alpha q^{k+n}}+\frac{\beta q^k}{1-\beta q^{k+n}}\right)
\nonumber \\
&\quad +(1-q^n)^2 (1-\gamma q^{n+k})\frac{\alpha\beta q^{2k}}{(1-\alpha q^{k+n})(1-\beta q^{k+n})}.
\end{align}
Unfortunately the right hand side of \eqref{RR2} does not give $\tilde{H}_{k+1}$, 
because of the last term, which vanishes only when $\alpha=0$ or $\beta=0$, {\it i.e.} nothing but the Ramanujan case. 

In conclusion, it is unable to extend the Ramanujan type continued fractions with parameters 
from two $(\beta, \gamma)$ to three $(\alpha,\beta,\gamma)$. \\

\noindent
{\bf (Continued fractions for two parameters $(\beta, \gamma)$)}\quad
Let us consider the analogous derivation of continued fraction in \S 2.2, of Gauss-Heine type, 
by using (0) and $(3)'$, 
\begin{align}
&(0)\quad [\beta, \gamma; x]_k-[\beta,\gamma;x]_{k+1}
\nonumber \\
&\qquad\qquad=q^{k+1} x\ \frac{1-\beta}{1-\gamma}\cdot [\beta q,\gamma q;x]_{k+1},
\\
&(3)'\quad 
[\beta q,\gamma q;x]_k-[\beta,\gamma; x]_k
\nonumber \\
&\qquad\qquad=q^{k+1}x\ \frac{(\beta-\gamma)}{(1-\gamma)(1-\gamma q)}\cdot[\beta q,\gamma q^2; x]_{k+1}.
\end{align}

It should be noted $[\beta,\gamma;x]_0=F_0,\ [\beta,\gamma;x]_1=H_0$, and
\begin{align}
&\frac{H_0}{F_0}=\frac{[\beta,\gamma;x]_1}{[\beta,\gamma;x]_0}
=\frac{e_0}{1-}\ 
\frac{e_1x}{\frac{[\beta, \gamma; x]_1}{[\beta q,\gamma q; x]_1}},
\\
&\qquad\text{with}\quad e_0=1,\ e_1=-\frac{q(1-\beta)}{(1-\gamma)}.
\nonumber
\end{align}

Next we have by using $(3)'$ with $k=1$,
\begin{align}
&\frac{[\beta q, \gamma q; x]_1}{[\beta,\gamma; x]_1}=
\frac{1}{1-}\ \frac{e_2 x}{\frac{[\beta q, \gamma q; x]_1}{[\beta q,\gamma q^2; x]_2}}
\\
&\qquad\text{with}\quad e_2=-\frac{q^2(\gamma-\beta)}{(1-\gamma)(1-\gamma q)}.
\nonumber
\end{align}
Such procedure can be repeated infinitely, and we obtain the same result of \eqref{RMC}, which was derived by Muir formula. 


\section{Examples of triplicity and applications}
\setcounter{equation}{0}


\subsection{Triplicity of Ramanujan type identities}

Andrews gave totally five examples of Ramanujan type continued fractions in \cite{Andrews2} from (7.10) to (7.14), 
It should be mentioned that he gave other different types of examples also. 
General method of making Ramanujan type continued fractions was already discussed in \S 2.4. 
These five concrete examples are themselves interesting to be examined, especially concerning to 
the triplicity property.   
Among them three have triplicity relations, and other two lack infinite product part $P$ unfortunately. \\

\noindent
{\bf (Example 1)}\quad The first one, (7.10) of Andrews, is given by
\begin{align}
&\prod_{n=0}^\infty \frac{(1-q^{4n+1})(1-q^{4n+3})}{(1-q^{4n+2})^2}
=\frac{\sum_{n=0}^\infty \frac{q^{n(n+1)}}{(q^2; q^2)_n}}{\sum_{n=0}^\infty \frac{q^{n^2}}{(q^2; q^2)_n}}
\nonumber \\
&\qquad =\frac{1}{1+}\ \frac{q}{1+}\ \frac{q^2+q}{1+}\ \frac{q^3}{1+}\ \cdots,\qquad \text{(R)}
\label{ex1_1}\\
&\qquad =\frac{1}{1+}\ \frac{q}{1+q+}\ \frac{q^2}{1+q^2+}\ \frac{q^3}{1+q^3+}\ \cdots,\qquad \text{(M)}
\label{ex1_2}
\end{align}
according to Ramanujan formula \eqref{RamanujanF} ($a=0,\ b=1,\ \lambda=1$), and Muir formula \eqref{MuirF}.  
We have added \eqref{ex1_2} to show non-uniqueness of $C$, where (R) or (M) stands for the method by Ramanujan or Muir. 

In fact Muir formula \eqref{MuirF} gives \eqref{ex1_2} by setting $x=1$ in the identity
\begin{align}
\frac{1+\sum_{n=1}^\infty \frac{q^{n(n+1)}}{(q^2; q^2)_n}\ x^n}{1+\sum_{n=1}^\infty \frac{q^{n^2}}{(q^2; q^2)_n}\ x^n}
=\frac{e_0}{1-}\ \frac{e_1x}{1-}\ \frac{e_2x}{1-}\ \frac{e_3x}{1-}\ \frac{e_4x}{1-}\ \cdots,
\end{align}
whose coefficients are $e_0=1$ and 
\begin{align}
&e_1=-\frac{q}{1+q},\ e_2=-\frac{q^2}{(1+q)(1+q^2)},\ e_3=-\frac{q^3}{(1+q^2)(1+q^3)},\ \cdots, 
\nonumber \\
&e_n=-\frac{q^n}{(1+q^{n-1})(1+q^n)},\quad (n \geq 2)
\end{align}
which is rewritten as \eqref{ex1_2} in the standard expression. \\

\noindent 
{\bf (Example 2)}\quad The next one, (7.11) of Andrews, is the case of Ramanujan formula \eqref{RamanujanF} 
($a=0, b=-1, \lambda=1$), and is given by
\begin{align}
\frac{\sum_{n=0}^\infty \frac{q^{n(n+1)}}{(q; q)_n^2}}{\sum_{n=0}^\infty \frac{q^{n^2}}{(q; q)_n^2}}
&=\frac{1}{1+}\ \frac{q}{1+}\ \frac{q^2-q}{1+}\ \frac{q^3}{1+}\ \frac{q^4-q^2}{1+}\ \cdots, \qquad \text{(R)}
\label{ex2_1}
\\
&=\frac{1}{1+}\ \frac{q}{1-q+}\ \frac{q^2}{1-q^2+}\ \frac{q^3}{1-q^3+}\ \cdots,\qquad \text{(M)}
\label{ex2_2}
\end{align}
where the last line is added, because Muir formula \eqref{MuirF} gives 
\begin{align}
&\frac{\sum_{n=0}^\infty \frac{q^{n(n+1)}}{(q; q)_n^2}x^n}{\sum_{n=0}^\infty \frac{q^{n^2}}{(q; q)_n^2}x^n}
=\frac{e_0}{1-}\ \frac{e_1x}{1-}\ \frac{e_2x}{1-}\ \frac{e_3x}{1-}\ \cdots,
\nonumber \\
&e_0=1,\ e_1=-\frac{q}{1-q},\ e_2=-\frac{q^2}{(1-q)(1-q^2)},\ \cdots,
\nonumber \\
&e_n=-\frac{q^n}{(1-q^{n-1})(1-q^n)}, \quad (n\geq 2).
\end{align}
Equation \eqref{ex2_2} is the standard expression of this result.

The continued fraction in the right hand side of \eqref{ex2_1} is a well known result. It is $x=-1$ case of 
\begin{align}
\sum_{n=0}^\infty \ q^{n(n+1)/2}\ x^n 
=\frac{1}{1-}\ \frac{qx}{1-}\ \frac{(q^2-q)x}{1-}\ \frac{q^3x}{1-}\ \frac{(q^4-q^2)x}{1-}\ \cdots,
\label{Eisenstein}
\end{align}
which implies the series part $S$ of this example is also not unique. 
Identity \eqref{Eisenstein} was found by Eisenstein in 1845 \cite{Eisenstein}. 
Gauss probably knew this also, because Entry 58 in 1797 of his diary is the case of $x=-1$. 

The infinite product expression is not known unfortunately except the case of $x=1$,
\begin{align}
\prod_{n=1}^\infty \left( \frac{1-q^{2n}}{1-q^{2n-1}} \right)=\sum_{n=0}^\infty q^{n(n+1)/2}
=\frac{1}{1-}\ \frac{q}{1-}\ \frac{q^2-q}{1-}\ \frac{q^3}{1-}\ \cdots,
\end{align}
which is the case of $\beta=q$ of another interesting identity of $P=S$ type
\begin{align}
\prod_{n=1}^\infty \left( \frac{1-\beta q^{2n-1}}{1-q^{2n-1}}\right)=
\sum_{n=0}^\infty q^{n(n+1)/2}\cdot\frac{(\beta; q)_n}{(q; q)_n},
\end{align}
the right hand side of which is, by the way, also the case $\gamma=0, x=1$ of $F_0(\beta,\gamma; q, x)$ 
introduced in \S 2.4, or the case $x=1$ of ${}_1\Phi_0(\beta; q, x)$ introduced in \S 2.5. 
Further the left hand side can also be written as
\begin{align}
\prod_{n=0}^\infty \left( \frac{1-\beta q\cdot (q^2)^{n}}{1-q\cdot (q^2)^{n}}\right)=
\sum_{n=0}^\infty q^n \frac{(\beta; q^2)_n}{(q^2;q^2)_n}
\end{align}
by use of Cauchy's identity \eqref{Cauchy}. This implies therefore 
\begin{align}
\sum_{n=0}^\infty q^{n(n+1)/2}\frac{(\beta;q)_n}{(q;q)_n}=\sum_{n=0}^\infty q^n\frac{(\beta; q^2)_n}{(q^2;q^2)_n}
\end{align}
which is related to a special case of Rogers-Fine identity discussed later. \\

\noindent
{\bf (Example 3)}\quad The third one, (7.12) of Andrews, is the case of Ramanujan formula \eqref{RamanujanF} 
($q\rightarrow q^2, a=q^{-1}, b=1, \lambda=1$), given by
\begin{align}
&\prod_{n=0}^\infty\frac{(1-q^{6n+1})(1-q^{6n+5})}{(1-q^{6n+3})^2}
=\frac{1+\sum_{n=1}^\infty \frac{q^{n(n+2)} (-q;q^2)_n}{(q^4; q^4)_n}}
{1+\sum_{n=1}^\infty \frac{q^{n^2} (-q; q^2)_n}{(q^4; q^4)_n}}
\nonumber \\
&\qquad\qquad =\frac{1}{1+}\ \frac{q+q^2}{1+}\ \frac{q^2+q^4}{1+}\ \frac{q^3+q^6}{1+}\ \cdots,\qquad \text{(R)}
\\
&\qquad\qquad =\frac{1}{1+}\ \frac{q(1+q)}{1+q^2+}\ \frac{q^4(1-q)}{1+q^4+}\ \frac{q^3(1+q^3)}{1+q^6+}\ \cdots.
\qquad \text{(M)}
\label{A7_12_M}
\end{align}
Muir formula \eqref{MuirF} gives
\begin{align}
&e_0=1,\ e_1=-\frac{q(1+q)}{1+q^2},\nonumber \\
&e_2=-\frac{q^4(1-q)}{(1+q^2)(1+q^4)},\ 
e_3=-\frac{q^3(1+q^3)}{(1+q^4)(1+q^6)},\ \cdots,
\nonumber
\end{align}
whose general expression is given by 
\begin{align}
\begin{split}
&e_{2m}=-\frac{q^{4m}(1-q^{2m-1})}{(1+q^{4m-2})(1+q^{4m})},\\ 
&e_{2m+1}=-\frac{q^{2m+1}(1+q^{2m+1})}{(1+q^{4m})(1+q^{4m+2})}.
\end{split}
\end{align}
for $m\geq 1$, which gives \eqref{A7_12_M} in the standard expression. \\

\noindent
{\bf (Example 4)}\quad The fourth one, (7.13) of Andrews, is the case of Ramanujan formula 
($q\rightarrow q^2, a=q^{-1}, b=0, \lambda=1)$, given by
\begin{align}
&\prod_{n=0}^\infty\frac{(1-q^{8n+1})(1-q^{8n+7})}{(1-q^{8n+3})(1-q^{8n+5})}
=\frac{1+\sum_{n=1}^\infty \frac{q^{n(n+2)} (-q;q^2)_n}{(q^2; q^2)_n}}
{1+\sum_{n=1}^\infty \frac{q^{n^2} (-q; q^2)_n}{(q^2; q^2)_n}}
\nonumber \\
&\qquad\qquad =\frac{1}{1+}\ \frac{q^2+q}{1+}\ \frac{q^4}{1+}\ \frac{q^6+q^3}{1+}\ \cdots,\quad(\text{R,\ M})
\end{align}
where the coefficients are generally 
\begin{align}
e_{2m}=-q^{4m},\quad e_{2m+1}=-(q^{4m+2}+q^{2m+1}),\quad (m\geq 1),
\end{align}
which are the same as derived by Muir formula \eqref{MuirF}, because Ramanujan of $b=0$ is equivalent to 
Muir.\\

\noindent
{\bf (Example 5)}\quad The last one, (7.14) of Andrews, is the case of Ramanujan formula 
($q\rightarrow q^2, a=-q^{-1},\ b=-1,\ \lambda=1$), given by
\begin{align}
&\frac{1+\sum_{n=1}^\infty (-1)^n \frac{q^{n(n+2)}\ (q; q^2)_n}{(q^2; q^2)_n^2}}
{1+\sum_{n=1}^\infty (-1)^n \frac{q^{n^2}\ (q; q^2)_n}{(q^2; q^2)_n^2}}
\nonumber \\
&\qquad\qquad =\frac{1}{1+}\ \frac{q^2-q}{1+}\ \frac{q^4-q^2}{1+}\ \frac{q^6-q^3}{1+}\ \cdots.
\qquad \text{(R)}
\label{A7_14}
\end{align}
According to Andrews \cite{Andrews2}, the left hand side (LHS) is known to be rewritten also by a single series such as
\begin{align}
\text{LHS}=\sum_{n=0}^\infty (-1)^n\frac{q^{n(n+1)}}{(q; q^2)_{n+1}}=\frac{1}{1-q}-\frac{q^2}{(1-q)(1-q^3)}+\cdots, 
\label{7_14}
\end{align}
which implies non-uniqueness of series part $S$. The product part $P$ is not known unfortunately. 

If we apply Muir-Rogers formula \eqref{MuirRogersF} to \eqref{7_14}, we have
\begin{align}
&\text{LHS}=\frac{e_0}{1-}\ \frac{e_1}{1-}\ \frac{e_2}{1-}\ \cdots,\qquad \text{(M)}
\nonumber \\
&\qquad e_0=\frac{1}{1-q},\ e_1=-\frac{q^2}{1-q^3},\ e_2=\frac{q^2(1-q^2)}{(1-q^3)(1-q^5)},\ \cdots,
\nonumber
\end{align}
where the coefficients are in general given by ($m\geq 1$)
\begin{align}
\begin{split}
e_{2m}&=+\frac{q^{2m}(1-q^{2m})}{(1-q^{4m-1})(1-q^{4m+1})},\\
e_{2m+1}&=-\frac{q^{4m+2}(1-q^{2m+1})}{(1-q^{4m+1})(1-q^{4m+3})},
\end{split}
\end{align}
which certainly disagree with those of \eqref{A7_14}, because the latter is derived by Ramanujan formula. 

On the contrary we can find another series which has the same coefficients of \eqref{A7_14} by using 
{\it inversion procedure} for Muir-Rogers formula, explained in \S 2.3. From $e_0=1,\ e_1=q-q^2,\ e_2=q^2-q^4,\ \cdots$, we have   
$\alpha_0=1,\ \alpha_1=q-q^2,\ \alpha_2=(q-q^2)(q^2-q^4),\ \cdots$, and finally 
we obtain $c_0=1,\ c_1=q(1-q),\ c_2=q^2(1-q)(1-q^3),\ \cdots$. 
The result is given by the identity
\begin{align}
\sum_{n=0}^\infty q^n(q; q^2)_n x^n=
\frac{1}{1+}\ \frac{(q^2-q)x}{1+}\ \frac{(q^4-q^2)x}{1+}\ \frac{(q^6-q^3)x}{1+}\ \cdots.
\end{align}
Therefore in conclusion  we have curious identity that three series exist for the common single continued fraction,
\begin{align}
\frac{1}{1+}\ \frac{q^2-q}{1+}\ \frac{q^4-q^2}{1+}\ \cdots
&=\frac{1+\sum_{n=1}^\infty (-1)^n \frac{q^{n(n+2)}\ (q; q^2)_n}{(q^2; q^2)_n^2}}
{1+\sum_{n=1}^\infty (-1)^n \frac{q^{n^2}\ (q; q^2)_n}{(q^2; q^2)_n^2}}
\nonumber \\
&=\sum_{n=0}^\infty (-1)^n\frac{q^{n(n+1)}}{(q; q^2)_{n+1}}
\nonumber \\
&=\sum_{n=0}^\infty q^n(q; q^2)_n,
\label{ex5_last}
\end{align}
which shows explicitly non-uniqueness of the series part $S$ for a given continued fraction $C$. 

\subsection{Answers to Gauss problems}

As was mentioned in \S 1 before, Gauss raised two problems in Entry 7 of his diary \cite{GaussD}. 
In this subsection we wish to give answers to them. \\

\noindent
{\bf(Gauss problem 1)}\quad Find the value of series 
\begin{align}
\sum_{n=0}^\infty (-1)^n q^{n(n+1)/2}=1-q+q^3-q^6+q^{10}-\cdots,
\label{GaussP_1}
\end{align}
for $q=2$ case. 

To find this value, we use the Rogers-Fine identity \cite{RogersFine1, RogersFine2}
\begin{align}
\sum_{n=0}^\infty x^n q^{n(n+1)/2} = \sum_{n=0}^\infty \frac{(xq; q^2)_n}{(xq^2; q^2)_n}\cdot (xq)^n,
\label{RogersFine}
\end{align}
at $x=-1$, which becomes
\begin{align}
\sum_{n=0}^\infty (-1)^n q^{n(n+1)/2}=1+\sum_{n=1}^\infty (-1)^n q^n\cdot
\frac{(1+q)\cdots(1+q^{2n-1})}{(1+q^2)\cdots(1+q^{2n})}.
\end{align}
Obviously both sides are convergent for $|q|<1$. But our problem is to discuss the case of $q>1$, 
although the left hand side is divergent so is {\it not defined} for $q>1$.  
However by substituting $q=p^{-1}$ in the right hand side, we have an equality
\begin{align}
\sum_{n=0}^\infty (-1)^n q^{n(n+1)/2}=
1+\sum_{n=1}^\infty (-1)^n \frac{(1+p)\cdots(1+p^{2n-1})}{(1+p^2)\cdots(1+p^{2n})},
\label{Gauss1_equality}
\end{align}
the right hand side of which is {\it almost convergent} when $p<1$, 
because it is similar to the alternating series of $Q_\infty$
\begin{align}
\sum_{n=0}^\infty (-1)^n Q_\infty,\qquad Q_\infty=\prod_{n=1}^\infty \left(\frac{1+p^{2n-1}}{1+p^{2n}}\right),
\label{Qinfinity}
\end{align}
whose sum is given by $\frac{1}{2}Q_\infty$ in the sense of Cez\`{a}ro sum \cite{HardyD}.
Therefore by using \eqref{Gauss1_equality} we can {\it define and replace} the left hand side for $q>1$ 
by the right hand side for $p=q^{-1}<1$. 

Let us take a short break here. 
The logic in the above treatment is analogous to the following situation. It is well known that
\begin{align}
1+x+x^2+\cdots=\frac{1}{1-x},
\nonumber
\end{align}
for $|x|<1$. For $|x|>1$ the left hand side series is divergent and {\it not defined}, 
and so let us {\it define} it by the same function such that
\begin{align}
\frac{1}{1-x}=-\left(\frac{1}{x}+\frac{1}{x^2}+\frac{1}{x^3}+\cdots\right),
\nonumber
\end{align}
which can be convinced by Euler's equality (see p.253 of \cite{Bourbaki})
\begin{align}
\sum_{n=-\infty}^\infty x^n=0,
\end{align}
which is understood in modern terminology of {\it distributions} (L. Schwartz) or {\it hyper-functions} (M. Sato) 
by
\begin{align}
&\sum_{n=-\infty}^\infty x^n=0\quad (x\neq 1),\qquad =\infty\quad (x=1)
\nonumber \\
&\Longleftrightarrow\quad 
\sum_{n=-\infty}^\infty e^{2\pi i n y}=\sum_{m=-\infty}^\infty \delta(y-m),
\end{align}
where $\delta(y-m)$ is Dirac's delta function, the latter equality is nothing but a generator of the Poisson summation formula.

Now let us return to our own problem. Result of numerical experiments of finite sum for $p=0.5\ (q=2)$ is summarized as follows.
\begin{align}
S_N&= 1+\sum_{n=1}^N (-1)^n \frac{(1+p)\cdots(1+p^{2n-1})}{(1+p^2)\cdots(1+p^{2n})}
\\
\begin{split}
&= +1.0759457568\cdots\quad(N=100),\\
&= -0.2208954963\cdots\quad(N=101).
\end{split}
\end{align}
Since $Q_\infty\fallingdotseq S_{N}-S_{N+1}$ for large even $N$, we can check the equality 
\begin{align}
Q_\infty=1.296841253\cdots=1.0759457568\cdots+0.2208954963\cdots
\end{align}
by computing $Q_\infty$ independently by \eqref{Qinfinity}. 
By the way, the above values are convergent already for $N=50$ up to the given ten digits. 

In other words, we can say for large even $N$ 
\begin{align}
S_N\fallingdotseq Lim+\frac{1}{2} Q_\infty,\quad
S_{N+1}\fallingdotseq Lim-\frac{1}{2} Q_\infty,
\end{align}
and the value of $Lim$ is given by the average 
\begin{align}
Lim\fallingdotseq\frac{1}{2}\left(S_N+S_{N+1}\right)=0.4275251302\cdots,
\end{align}
for $p=0.5,\ N=100$. 

In conclusion the answer to (Gauss problem 1) is given by the above
\begin{align}
1-2+2^3-2^6+\cdots=0.4275251302\cdots,
\end{align}
which is our result. \\

Some remarks should be added here. By the same manner,
setting $q\rightarrow q^2$ and $x=-q^{-1}$ in \eqref{RogersFine}, we can derive the equality
\begin{align}
\sum_{n=0}^\infty (-1)^n q^{n^2}&=1+\sum_{n=1}^\infty (-1)^n q^n \frac{(1+q)\cdots(1+q^{4n-3})}{(1+q^3)\cdots(1+q^{4n-1})}
\nonumber \\
&=1+\sum_{n=1}^\infty (-1)^n p^n \frac{(1+p)\cdots(1+p^{4n-3})}{(1+p^3)\cdots(1+p^{4n-1})},
\label{duality}
\end{align}
where $p=q^{-1}$ as usual. It is interesting that two right hand sides have the same form, 
and therefore two series of $q<1$ and $q>1$ have the same value as their sum, that is, 
this example has a kind of {\it self-dual} property. In fact for $q=0.5\ (p=2)$ and $p=0.5\ (q=2)$, all take the same value 
numerically
\begin{align}
&1-2^{-1}+2^{-4}-2^{-9}+\cdots=1-2+2^4-2^9+\cdots
\nonumber \\
&\qquad\qquad\qquad\quad =0.5605621040 \cdots
\end{align}
which is computed for a large enough finite sum $N=100$. 
Here it must be stressed again that we regard the case of $q>1$ is {\it defined} by the $p$ series in the right hand side of 
\eqref{duality}.

We can obtain similar identities for every polygonal number series. For example the case of pentagonal numbers is given by
\begin{align}
&\sum_{n=0}^\infty (-1)^n q^{n(3n-1)/2}=1+\sum_{n=1}^\infty (-1)^n q^n\frac{(1+q)\cdots(1+q^{6n-5})}
{(1+q^4)\cdots(1+q^{6n-2})}
\nonumber \\
&\qquad =1+\sum_{n=1}^\infty (-1)^n p^{2n} \frac{(1+p)\cdots(1+p^{6n-5})}
{(1+p^4)\cdots(1+p^{6n-2})},\quad (p=q^{-1})
\end{align}
which implies two series in the right hand sides for $q<1$ and $p<1$ have convergent values. 
For example the values for $q=p=0.5$ are computed as follows
\begin{align}
\begin{split}
&1-2^{-1}+2^{-5}-2^{-12}+\cdots=0.5310060977\cdots,\\
&1-2+2^5-2^{12}+\cdots=0.7181272344\cdots.
\end{split}
\end{align}
The latter is our result of divergent series for $q=2$, computed by the right hand side at $p=0.5$. \\

\noindent
{\bf(Gauss problem 2)}\quad Find the value of series 
\begin{align}
\sum_{n=0}^\infty (q;q)_n=1+(1-q)+(1-q)(1-q^2)+\cdots.
\label{GaussP_2}
\end{align}
for $q=2$ case.

Before answering this question, let us discuss about other related problems first. 
It is easy to notice that \eqref{GaussP_2} is divergent even if $|q|<1$,  
because the value of $(q; q)_\infty=\prod_{n=1}^\infty(1-q^n)\equiv Q_\infty$ 
is non-zero finite except the case of $q=1$, or the radicals to $q^N=1$ for complex $q$ and integer $N$.
We hope the same symbol $Q_\infty$ does not make confusion with the previous problem.

On the other hand, the alternating series
\begin{align}
\sum_{n=0}^\infty (-1)^n (q;q)_n=1-(1-q)+(1-q)(1-q^2)-\cdots
\end{align}
can be regarded convergent in the same sense of the previous problem, because it is again similar to 
$\sum_{n=0}^\infty (-1)^n Q_\infty$, and the value of $Lim$ can be computed by the average as before. 
In fact we have  
\begin{align}
\begin{split}
&S_N=\sum_{n=0}^N (-1)^n (q;q)_n,
\\
&Lim\fallingdotseq\frac{1}{2}\left(S_{N}+S_{N+1}\right)=0.7039282729\cdots,
\end{split}
\end{align} 
for $q=0.5$ and $N=100$ finite sum. 

The same result can be obtained by using Heine's transformation identity
\begin{align}
&\sum_{n=0}^\infty \frac{(\alpha; q)_n(\beta;q)_n}{(\gamma;q)_n(q;q)_n}\ x^n
\nonumber \\
&\qquad=
\frac{(\gamma/\beta;q)_\infty(\beta x;q)_\infty}{(\gamma;q)_\infty(x;q)_\infty}
\sum_{n=0}^\infty \frac{(\beta;q)_n(\alpha\beta x/\gamma;q)_n}{(\beta x;q)_n(q;q)_n}\left(\frac{\gamma}{\beta}\right)^n,
\label{HeineTr}
\end{align}
which is (1.4.5) in \cite{GasperRahman} by Gasper-Rahman.
If we set $\alpha=\beta=q, \gamma=0$, we have
\begin{align}
\sum_{n=0}^\infty (q; q)_n x^n=
\frac{1}{1-x}\left(\sum_{n=0}^\infty (-1)^n \frac{q^{n(n+1)/2}}{(xq; q)_n}\ x^n\right),
\label{GaussP_3}
\end{align}
which shows clearly that the left hand side diverges when $x=1$. 
On the contrary, if we set $x=-1$  we obtain the equality
\begin{align}
\sum_{n=0}^\infty (-1)^n (q;q)_n=\frac{1}{2}\left(1+\sum_{n=1}^\infty \frac{q^{n(n+1)/2}}{(1+q)\cdots(1+q^n)}\right).
\label{GR}
\end{align}
The right hand side gives directly the same value $0.7039282729\cdots$ when $q=0.5$ without taking average. 
It is amusing to find \eqref{GR} gives directly
\begin{align}
\sum_{n=0}^\infty (-1)^n=1-1+1-1+\cdots=\frac{1}{2}
\label{Cezaro}
\end{align}
by substituting $q=0$, which was interpreted as a C\`{e}zaro sum before.

Let us consider next the case $q>1,\ p=q^{-1}<1$, and try to estimate the parentheses part in the right hand side of 
\eqref{GaussP_3} for $x=1$. It is given by 
\begin{align}
\left( \cdots \right)&=1+\sum_{n=1}^\infty (-1)^n \frac{q^{n(n+1)/2}}{(1-q)\cdots(1-q^n)}
\nonumber \\
&=1+\sum_{n=1}^\infty \frac{1}{(1-p)\cdots(1-p^n)},
\end{align}
the right hand side of which is divergent for $p<1$ because every term in the series is positive definite, implying that
the right hand side of \eqref{GaussP_3} is {\it doubly divergent} like $\infty/0=\infty^2$ for $x=1$. 

In conclusion, our answer to (Gauss problem 2) is unexpectedly 
\begin{align}
\sum_{n=0}^\infty (2; 2)_n=1-1+1\cdot 3-1\cdot 3\cdot 7+\cdots=\text{divergent},
\end{align}
even though the middle part looks like an alternating series. \\

It is desirable to find an alternative problem which gives a finite value to similar divergent series. 
The next is such a problem. \\

\noindent
{\bf (Gauss problem 3)}\quad Find the value of
\begin{align}
\sum_{n=0}^\infty q^n(q;q^2)_n=1+q(1-q)+q^2(1-q)(1-q^3)+\cdots
\end{align}
for $q=2$ case. Let us call this also Gauss problem, although he did not submit this.

Obviously this series is convergent for $|q|<1$, and therefore we focus on the case of $q>1$. 
According to \eqref{ex5_last} of (Example 5) in \S 3.1, we have equalities
\begin{align}
\sum_{n=0}^\infty q^n(q;q^2)_n&=\sum_{n=0}^\infty (-1)^n\frac{q^{n(n+1)}}{(q; q^2)_{n+1}}
\nonumber \\
&=-\left(\sum_{n=0}^\infty \frac{p^{n+1}}{(p; p^2)_{n+1}}\right),\qquad (p=q^{-1}).
\label{GaussP3}
\end{align}
It is noticeable that the right hand side of the last line has over all negative sign, 
and it is certainly convergent for $|p|<1$. 
If we set $q=2,\ p=0.5$, the above equality \eqref{GaussP3} becomes
\begin{align}
&\sum_{n=0}^\infty 2^n(2;4)_n=1-2+4\cdot7-8\cdot 7\cdot 31+\cdots
\nonumber \\
&\qquad =-\left(\sum_{n=0}^\infty \frac{2^{-(n+1)}}{(2^{-1};2^{-2})_{n+1}}\right)=-2.1639450388\cdots,
\end{align}
which is our answer to (Gauss problem 3). \\

The numbers obtained here are all computed numerically and unfortunately give no information about their arithmetical 
meanings. Anyway our results in this subsection all remind us the famous equalities
\begin{align}
\begin{split}
&1+2+3+\cdots=\zeta(-1)=-\frac{1}{12},\quad \text{hence}\\
&1-2+3-\cdots=(1-2^2) \cdot\zeta(-1)=+\frac{1}{4},
\end{split}
\end{align}
which are typical examples of divergent series taking a finite value.  
Further we can derive 
\begin{align}
\begin{split}
&1+1+1+\cdots=\zeta(0)=-\frac{1}{2},\quad \text{because}\\
&\frac{1}{2}=1-1+1-1+\cdots=(1-2^1)\cdot\zeta(0),
\end{split}
\end{align}
the left hand side of which is nothing but \eqref{Cezaro}. 
Such values $\zeta(-1), \zeta(0)$ are usually derived by use of functional equation of zeta function
\begin{align}
\zeta(1-s)=\frac{2\cos(\pi s/2)}{(2\pi)^s}\ \Gamma(s)\zeta(s),
\end{align}
which was proved by Riemann (1859) \cite{Riemann} using {\it analytical continuation}.

In the above we used naively the equality
\begin{align}
\sum_{n=1}^\infty \frac{(-1)^{n-1}}{n^s}=(1-2^{1-s})\cdot\zeta(s),
\end{align}
which is the Fermi integral $F_s$, a counterpart of the Bose integral $B_s$
\begin{align}
\begin{split}
&B_s=\int_0^\infty \frac{x^{s-1}dx}{e^x-1}=\zeta(s)\ \Gamma(s),\\
&F_s=\int_0^\infty \frac{x^{s-1}dx}{e^x+1}=(1-2^{1-s})\cdot \zeta(s)\ \Gamma(s).
\end{split}
\end{align}
The latter can also be used for analytical continuation in the same manner as Riemann used the former.

The mechanism of computing divergent series in this subsection is based on {\it non-uniqueness} 
of series part $S$ in the triplicity. However it can be regarded also as {\it a kind of analytical continuation}, because 
one expression is considered as analytically continued the other one.

\subsection{Triplicity applied to divergent series}

In this subsection, five examples given in \S 3.1, related to Ramanujan's continued fractions, 
are re-examined with respect to sums of divergent series. 
Our idea to compute divergent series is to use {\it non-uniqueness} of series part $S$ in the triplicity $P=S=C$. 
Such examples are already given in the previous subsection to answer Gauss problems. 
For the concreteness, let us state our strategy explicitly. \\

{\it Suppose the series part is not unique: $S_1(q)=S_2(q)$, and $S_1(q)$ is convergent (divergent) for $q<1$ ($q>1$).
Then further if $S_2(q=p^{-1})\equiv \tilde{S}_2(p)$ is convergent for $p<1$, we can replace the value of divergent series $S_1(q)$ 
by the sum of convergent series $\tilde{S}_2(p)$. 
On the contrary, if $\tilde{S}_2(p)$ is still divergent, we conclude $S_1(q)$ is divergent for $q>1$.} \\

One comment is in order before entering each topics. Three examples among five have product part $P$ with a common property,
\begin{align}
P(q)=\prod \frac{(1-q^a)(1-q^b)}{(1-q^c)(1-q^d)},\qquad(a+b=c+d).
\end{align}
Since we notice, dividing numerator and denominator by $q^{a+b}=q^{c+d}$, 
\begin{align}
\frac{(1-q^a)(1-q^b)}{(1-q^c)(1-q^d)}&=\frac{(q^{-a}-1)(q^{-b}-1)}{(q^{-c}-1)(q^{-d}-1)}
=\frac{(1-q^{-a})(1-q^{-b})}{(1-q^{-c})(1-q^{-d})}
\nonumber \\
&=\frac{(1-p^a)(1-p^b)}{(1-p^c)(1-p^d)},\qquad(p=q^{-1})
\nonumber
\end{align}
we obtain finally 
\begin{align}
P(q)=\prod \frac{(1-q^a)(1-q^b)}{(1-q^c)(1-q^d)}=\prod \frac{(1-p^a)(1-p^b)}{(1-p^c)(1-p^d)})=P(p),
\label{P=P}
\end{align}
for $p=q^{-1}$. Therefore it seems obvious that we have {\it self-duality} such that $P(q)=P(q^{-1})$. 

However the situations are not so trivial, because the series part $S$ must be examined carefully, 
and we will find a contradiction against above expectation. Suppose that $S(q)$ is corresponding series part
\begin{align}
P(q)=\prod \frac{(1-q^a)(1-q^b)}{(1-q^c)(1-q^d)}=S(q),\qquad (q<1)
\nonumber
\end{align}
and when $q>1$ we define $\tilde{S}(p)$ by {\it substitution} such as
\begin{align}
\tilde{S}(p)=S(q=p^{-1}),\qquad (p<1)
\nonumber
\end{align}
then we encounter usually
\begin{align}
\tilde{S}(p)\neq S(q),\qquad (p=q<1)
\end{align}
which implies $P(q^{-1})\neq P(q)$. 
Such phenomenon was seen already in the previous subsection, except \eqref{duality} that was indeed self-dual, 
and will be seen explicitly also in cases below. 
Therefore our formal rewriting \eqref{P=P} above contains some logical fallacy. \\

\noindent
{\bf (Case 1)}\quad This is the topic of Example 1 in \S 3.1, which corresponds to (7.10) of Andrews \cite{Andrews2},
and its triplicity is as follows (continued fraction part $C$ is omitted),
\begin{align}
\prod_{n=0}^\infty \frac{(1-q^{4n+1})(1-q^{4n+3})}{(1-q^{4n+2})^2}
=\frac{1+\sum_{n=1}^\infty \frac{q^{n(n+1)}}{(q^2; q^2)_n}}{1+\sum_{n=1}^\infty \frac{q^{n^2}}{(q^2; q^2)_n}}
\label{Case1_1}
\end{align}
where it should be stressed that the above formulas are for $|q|<1$. 
Now let us write 
\begin{align}
\begin{split}
&S(q)=\frac{1+\sum_{n=1}^\infty \frac{q^{n(n+1)}}{(1-q^2)\cdots(1-q^{2n})}}
{1+\sum_{n=1}^\infty \frac{q^{n^2}}{(1-q^2)\cdots(1-q^{2n})}},
\\
&\tilde{S}(p)=S(q=p^{-1})=
\frac{1+\sum_{n=1}^\infty (-1)^n \frac{1}{(1-p^2)\cdots(1-p^{2n})}}
{1+\sum_{n=1}^\infty (-1)^n \frac{p^n}{(1-p^2)\cdots(1-p^{2n})}}.
\end{split}
\end{align}
Two formulas should have the same value because the latter is simply a rewriting of the former when $q<1\ (p>1)$. 
And according to our strategy we {\it define} $S(q)\ (q>1)$ by $\tilde{S}(p)\ (p<1)$, which are all convergent 
except the numerator of $\tilde{S}(p)$, which is convergent in the sense of C\`{e}zaro.
 
Numerical computations for $q=0.5$ and $p=0.5$ for finite $N=100$ sums give respectively
\begin{align}
\begin{split}
&S(q)=0.7711044027\cdots,
\\
&\tilde{S}(p)=0.6484206265\cdots.
\end{split}
\end{align}
Therefore as was announced before, $S(q)\neq\tilde{S}(p)$ for $q=p<1$, although both are finite.

We conclude that this case does not have self-duality contrary to the simple expectation \eqref{P=P}. \\

\noindent
{\bf (Case 2)}\quad This is Example 2 of \S 3.1, which corresponds to (7.11) of Andrews \cite{Andrews2}. 
This case was used for (Gauss problem 1) and gave satisfactory result, which is not repeated here. 
For the sake of completeness, let us give only various compatible series part $S$. 

\begin{align}
&\sum_{n=0}^\infty (-1)^n q^{n(n+1)/2}
=1+\sum_{n=1}^\infty (-1)^n q^n  \frac{(1+q)\cdots(1+q^{2n-1})}{(1+q^2)\cdots(1+q^{2n})}
\nonumber \\
&\quad =1+\sum_{n=1}^\infty (-1)^n \frac{(1+p)\cdots(1+p^{2n-1})}{(1+p^2)\cdots(1+p^{2n})},
\quad (p=q^{-1})
\end{align}
which are convenient forms for $q<1\ (p>1)$. 
As we know from the above this case does not have self-duality, since the latter lacks $p^n$ coefficient, 
see also \eqref{duality}. Infinite product part $P$ is not known.\\

\noindent
{\bf (Case 3)}\quad This is Example 3 in \S 3.1, which is (7.12) of Andrews \cite{Andrews2}, 
whose triplicity is given by (continued fraction part $C$ is again omitted)
\begin{align}
\prod_{n=0}^\infty \frac{(1-q^{6n+1})(1-q^{6n+5})}{(1-q^{6n+3})^2}
=\frac{1+\sum_{n=1}^\infty q^{n(n+2)}\frac{(-q;q^2)_n}{(q^4;q^4)_n}}
{1+\sum_{n=1}^\infty q^{n^2} \frac{(-q;q^2)_n}{(q^4;q^4)_n}},
\end{align}
for $|q|<1$. Then let us write 
\begin{align}
\begin{split}
&S(q)=\frac{1+\sum_{n=1}^\infty q^{n(n+2)} \frac{(1+q)\cdots(1+q^{2n-1})}{(1-q^4)\cdots(1-q^{4n})}}
{1+\sum_{n=1}^\infty q^{n^2}\frac{(1+q)\cdots(1+q^{2n-1})}{(1-q^4)\cdots(1-q^{4n})}},
\\
&\tilde{S}(p)=S(q=p^{-1})=\frac{1+\sum_{n=1}^\infty(-1)^n \frac{(1+p)\cdots(1+p^{2n-1})}{(1-p^4)\cdots(1-p^{4n})}}
{1+\sum_{n=1}^\infty (-1)^n p^{2n}\frac{(1+p)\cdots(1+p^{2n-1})}{(1-p^4)\cdots(1-p^{4n})}}.
\end{split}
\end{align}
Similar to (case 1), above $S(q)$ are also convergent for $q<1$, and the numerator of $\tilde{S}(p)$ is C\`{e}zaro 
convergent for $p<1$. 

Numerical results of finite $N=100$ sums give respectively
\begin{align}
\begin{split}
&S(q)=0.6298180171\cdots,
\\
&\tilde{S}(p)=0.4072795451\cdots,
\end{split}
\end{align}
for $q=p=0.5$. As was mentioned, we have again $S(q)\neq\tilde{S}(p)$ for $q=p<1$, that is, our case is neither self-dual. \\

\noindent
{\bf (Case 4)}\quad This is Example 4 in \S 3.1, which is (7.13) of Andrews \cite{Andrews2}, 
whose triplicity is given by (continued fraction part $C$ is again omitted)
\begin{align}
\prod_{n=0}^\infty \frac{(1-q^{8n+1})(1-q^{8n+7})}{(1-q^{8n+3})(1-q^{8n+5})}
=\frac{1+\sum_{n=1}^\infty q^{n(n+2)}\frac{(-q;q^2)_n}{(q^2;q^2)_n}}
{1+\sum_{n=1}^\infty q^{n^2}\frac{(-q;q^2)_n}{(q^2;q^2)_n}},
\end{align}
for $|q|<1$ case. Let us write
\begin{align}
\begin{split}
&S(q)=\frac{1+\sum_{n=1}^\infty q^{n(n+2)}\frac{(1+q)\cdots(1+q^{2n-1})}{(1-q^2)\cdots(1-q^{2n})}}
{1+\sum_{n=1}^\infty q^{n^2}\frac{(1+q)\cdots(1+q^{2n-1})}{(1-q^2)\cdots(1-q^{2n})}}
\\
&\tilde{S}(p)=S(q=p^{-1})=\frac{1+\sum_{n=1}^\infty (-1)^n p^{-n(n+1)}\frac{(-p;p^2)_n}{(p^2;p^2)_n}}
{1+\sum_{n=1}^\infty (-1)^n p^{-n^(n-1)}\frac{(-p;p^2)_n}{(p^2;p^2)_n}}.
\end{split}
\end{align}
This $S(q)$ is convergent for $q<1$, however $\tilde{S}(p)$ for $p<1$ is {\it not convergent}, rather 
very rapidly diverges, which is analogous to (Gauss problem 2) in the previous subsection. 

For the concreteness let us give the numerical result of finite $N=100$ sum 
\begin{align}
\begin{split}
&S(q)=0.5844460945\cdots,\\
&\tilde{S}(p)=\text{divergent},
\end{split}
\end{align}
for $q=0.5$, and $p=0.5$. 
In conclusion we can say that this case also gives $S(q)\neq \tilde{S}(p)$, although  the latter is not finite. \\

\noindent
{\bf (Case 5)}\quad This is Example 5, which is (7.14) of Andrews, and was used for (Gauss problem 3) successfully. 
Details are referred to \S 3.2, and let us write only the compatible series below,
\begin{align}
&\frac{1+\sum_{n=1}^\infty q^{n(n+2)}\frac{(1-q)\cdots(1-q^{2n-1})}{(1-q^2)\cdots(1-q^{2n})}}
{1+\sum_{n=1}^\infty q^{n^2} \frac{(1-q)\cdots(1-q^{2n-1})}{(1-q^2)\cdots(1-q^{2n})}}
=\sum_{n=0}^\infty q^n (q; q^2)_n 
\nonumber \\
&\quad=\sum_{n=0}^\infty (-1)^n \frac{q^{n(n+1)}}{(q; q^2)_{n+1}}
=-\left(\sum_{n=0}^\infty \frac{p^{n+1}}{(p;p^2)_{n+1}}\right),\quad(p=q^{-1})
\end{align}
while the infinite product $P$ is not known unfortunately.\\

\noindent
{\bf (Relation to Jacobi's triple product identity)}\quad 
Changing the topic, we wish to discuss here another aspect of triplicity relation. 
The $P=S$ part in triplicity can be decoupled for some cases, for example the Rogers-Ramanujan identity \eqref{RRIdentity} has
\begin{align}
&\prod_{n=1}^\infty \frac{1}{(1-q^{5n-1})(1-q^{5n-4})}=1+\sum_{n=1}^\infty \frac{q^{n^2}}{(1-q)\cdots(1-q^n)},
\\
&\prod_{n=1}^\infty \frac{1}{(1-q^{5n-2})(1-q^{5n-3})}=1+\sum_{n=1}^\infty \frac{q^{n(n+1)}}{(1-q)\cdots(1-q^n)},
\end{align}
separately, which are (7.1.6) and (7.1.7) of \cite{Andrews}. 

Let us mention another facts, related to the above situation, which are given by the other identities
\begin{align}
&\prod_{n=1}^\infty(1-q^{5n})(1-q^{5n-1})(1-q^{5n-4})
\nonumber \\
&\qquad =1+\sum_{n=1}^\infty (-1)^n \left(q^{n(5n+3)/2}+q^{n(5n-3)/2}\right),
\label{Jacobi1}
\\
&\prod_{n=1}^\infty(1-q^{5n})(1-q^{5n-2})(1-q^{5n-3})
\nonumber \\
&\qquad =1+\sum_{n=1}^\infty (-1)^n \left(q^{n(5n+1)/2}+q^{n(5n-1)/2}\right),
\label{Jacobi2}
\end{align}
from which we have another expression for $S$ part of Rogers-Ramanujan identity such that
\begin{align}
&\prod_{n=1}^\infty\frac{(1-q^{5n-1})(1-q^{5n-4})}{(1-q^{5n-2})(1-q^{5n-3})}
=\frac{1+\sum_{n=1}^\infty \frac{q^{n(n+1)}}{(1-q)\cdots(1-q^n)}}
{1+\sum_{n=1}^\infty \frac{q^{n^2}}{(1-q)\cdots(1-q^n)}}
\nonumber \\
&\qquad =\frac{1+\sum_{n=1}^\infty (-1)^n\left(q^{n(5n+3)/2}+q^{n(5n-3)/2}\right)}
{1+\sum_{n=1}^\infty (-1)^n\left(q^{n(5n+1)/2}+q^{n(5n-1)/2}\right)}.
\end{align}
This equality can be tested numerically, and all values for $q=0.5$ coincide with $0.7099166943\cdots$, 
which are all convergent for finite $N=100$ product and sums.

Identities \eqref{Jacobi1} and \eqref{Jacobi2} are derived from the famous Jacobi's triple product identity, 
which is written by
\begin{align}
&\prod_{m=1}^\infty (1-p^{2m})(1+z^2p^{2m-1})(1+z^{-2}p^{2m-1})
\nonumber \\
&\qquad =
1+\sum_{m=1}^\infty p^{m^2}\left(z^{2m}+z^{-2m}\right).
\end{align}
By setting $p=q^{a/2},\ z^2=-q^{b/2}$ we have an identity
\begin{align}
&\prod_{m=1}^\infty (1-q^{am})(1-q^{am-(a-b)/2})(1-q^{am-(a+b)/2})
\nonumber \\
&\qquad =1+\sum_{m=1}^\infty (-1)^m\left(q^{m(am+b)/2}+q^{m(am-b)/2}\right).
\label{Elkies}
\end{align}
\eqref{Jacobi1} and \eqref{Jacobi2} are the case of $a=5,\ b=3$ and $a=5,\ b=1$ respectively. 

In the same manner by choosing appropriate values for $a$ and $b$, we can obtain another expression of $S$ part 
for the above cases (1), (3), and (4) which possess the same infinite product part $P$. 
Such values are respectively $a=4$ for (1), $a=6$ for (3), and $a=8$ for (4). 
The explicit formulas are omitted to write down here, because the author has no idea for what to use 
such relations, but see the last remark in the next section on $a=7$ case.


\section{Summary and concluding remarks}
\setcounter{equation}{0}

{\bf (Summary)}\quad The concept of triplicity is introduced and various examples are shown satisfying triplicity relation $P=S=C$, 
among infinite product $P$, infinite series $S$, and continued fractions $C$. 
Then its applications to computing sums of divergent series are given, including answers to historical Gauss problems. 

In this last section, we wish to give some remarks about questions left, new problems to be studied, and related topics 
to be examined. \\

\noindent
{\bf (Infinite dimensional determinants)}\quad 
Continued fraction $C$ in standard expression can be expressed by a quotient of two determinants such that
\begin{align}
&C=\frac{a_0}{b_0+}\ \frac{a_1}{b_1+}\ \frac{a_2}{b_2+}\ \cdots=a_0\cdot \frac{\Delta_1}{\Delta_0},
\label{Cdd}
\end{align}
where infinite dimensional determinants $\Delta_0$ and $\Delta_1$ are defined by 
\begin{align}
\begin{split}
&\Delta_0=\left|\begin{array}{ccccc}
b_0&a_1&0&0&\cdots\\
-1&b_1&a_2&0&\cdots\\
0&-1&b_2&a_3&\cdots\\
0&0&-1&b_3&\cdots\\
\cdots&\cdots&\cdots&\cdots&\cdots
\end{array}\right|,\ 
\Delta_1=\left|\begin{array}{ccccc}
b_1&a_2&0&0&\cdots\\
-1&b_2&a_3&0&\cdots\\
0&-1&b_3&a_4&\cdots\\
0&0&-1&b_4&\cdots\\
\cdots&\cdots&\cdots&\cdots&\cdots
\end{array}\right|,
\end{split}
\nonumber
\end{align}
which are generally expressed by
\begin{align}
\Delta_n=\left|\begin{array}{ccccc}
b_n&a_{n+1}&0&0&\cdots\\
-1&b_{n+1}&a_{n+2}&0&\cdots\\
0&-1&b_{n+2}&a_{n+3}&\cdots\\
0&0&-1&b_{n+3}&\cdots\\
\cdots&\cdots&\cdots&\cdots&\cdots
\end{array}\right|,\quad (n\geq 0).
\end{align}
It should be noted that they satisfy the iterative relations
\begin{align}
\begin{split}
&\Delta_0=b_0\Delta_1+a_1\Delta_2,\ \Delta_1=b_1\Delta_2+a_2\Delta_3,\ \cdots,\\
&\Delta_{n-1}=b_{n-1}\Delta_n+a_n\Delta_{n+1},\quad(n\geq 1)
\end{split}
\end{align}
from which we can reconstruct the original continued fraction $C$, {\it i.e.} the left hand side of \eqref{Cdd}.

The simplest example is the case of Rogers-Ramanujan identity. Let us introduce an infinite dimensional determinant
\begin{align}
\Delta(x)=\left|\begin{array}{ccccc}
1&qx&0&0&\cdots\\
-1&1&q^2x&0&\cdots\\
0&-1&1&q^3x&\cdots\\
0&0&-1&1&\cdots\\
\cdots&\cdots&\cdots&\cdots&\cdots
\end{array}\right|,
\end{align}
which is given by setting $a_0=b_0=1$, and $b_n=1,\ a_n=q^nx,\ (n\geq 1)$. 
It should be noted that $\Delta(x)=\Delta_0,\ \Delta(qx)=\Delta_1$, and $\Delta(0)=1$, 
therefore from \eqref{RRM} we have 
\begin{align}
\Delta(x)=F(x)=1+\sum_{n=1}^\infty \frac{q^{n^2}}{(q;q)_n} x^n,
\label{RRDet}
\end{align}
which implies that we have an expression of infinite dimensional determinant written by an explicit power series. 

Further we can say that all Ramanujan type series $G(a,\lambda;b,q)$ discussed in \S 2.4 also have the corresponding 
infinite dimensional determinant, because we have
\begin{align}
\frac{H_0}{F_0}=
\frac{a_0}{1+}\ \frac{a_1}{1+}\ \frac{a_2}{1+}\ \frac{a_3}{1-}\ \cdots,
\end{align}
where the coefficients are given by $a_0=1$ and
\begin{align}
&a_1=q(a+\lambda),\ a_2=q(b+\lambda q),\ a_3=q^2(a+\lambda q),\ \cdots,
\nonumber \\
&a_{2m}=q^m(b+\lambda q^m),\quad
a_{2m+1}=q^{m+1}(a+\lambda q^m),\ (m\geq 1).
\label{RC}
\end{align}
And therefore similarly we obtain
\begin{align}
F_0&=1+\sum_{n=1}^\infty \frac{q^{n(n+1)/2}}{(q;q)_n}\cdot
\frac{(a+\lambda)\cdots(a+\lambda q^{n-1})}{(1+bq)\cdots(1+bq^n)}
\nonumber \\
&=\left|\begin{array}{ccccc}
1&a_1&0&0&\cdots\\
-1&1&a_2&0&\cdots\\
0&-1&1&a_3&\cdots\\
0&0&-1&1&\cdots\\
\cdots&\cdots&\cdots&\cdots&\cdots
\end{array}\right|,
\label{RID}
\end{align}
where matrix elements are given by \eqref{RC}. Certainly as a special case of $a=0, \lambda=x, b=0$, we obtain \eqref{RRDet}.

Generally speaking, it is very difficult to rewrite forwardly (and backwardly also) from power series to 
infinite dimensional determinant like \eqref{RID}, such task Ramanujan achieved successfully.\\

\noindent
{\bf (Stieltjes integral)}\quad 
In 1894 Stieltjes wrote a unique paper on continued fractions \cite{Stieltjes}, where he introduced an integral now called with 
his name (see p.284 of \cite{Bourbaki}). 
In that paper he also raised and solved {\it the moment problem} which is related to the theories of probability and 
orthogonal polynomials \cite{Sego}. 

According to Van Vleck \cite{VanVleck}, theory of Stieltjes can be summarized as follows. 
Let us consider first a power series
\begin{align}
\begin{split}
&\frac{c_0}{z}-\frac{c_1}{z^2}+\frac{c_2}{z^3}-\frac{c_3}{z^4}+\cdots=\int_0^\infty \frac{d\Phi(u)}{z+u},\quad(z>0)
\\
&\quad\Longleftrightarrow\quad c_n=\int_0^\infty u^n d\Phi(u)
\end{split}
\end{align}
where $c_n$ is the $n$-th moment for probability density function $d\Phi(u)/du$. 
In the context of probability theory, the function $\Phi(u)$ must be non-decreasing positive function. 
For a better connection with the present paper's formulation  
we should have considered with the pivot variable $x=-z^{-1}$,
\begin{align}
\begin{split}
c_0+c_1x+c_2x^2+\cdots&=\int_0^\infty \frac{d\Phi(u)}{1-xu},\quad(x<0)
\\
&=\frac{e_0}{1-}\ \frac{e_1x}{1-}\ \frac{e_2x}{1-}\ \cdots,
\end{split}
\end{align}
which may remind us the triplicity of Wallis problem \eqref{WallisTriplicity}. 
In fact it is a typical example of Stieltjes integral, with $c_n=(a)_n$ and $d\Phi(u)/du=u^{a-1}e^{-u}/\Gamma(a)$. 
See also a note in \cite{RogersFine2} (p.32) on probabilistic interpretation of ${}_1\phi_1(aq, bq; q, t)$ as another example. 
The procedure to find $\Phi(u)$ from given $c_n$'s is called the moment problem, to which Stieltjes gave the answer.

Nextly Van Vleck \cite{VanVleck} introduces determinants $A_n, B_n$ 
defined by 
\begin{align}
\begin{split}
&A_n=\left|\begin{array}{cccc}
c_0&c_1&\cdots&c_{n-1}\\
c_1&c_2&\cdots&c_n\\
\cdots&\cdots&\cdots&\cdots\\
c_{n-1}&c_n&\cdots&c_{2n-2}\end{array}\right|,\\
&B_n=\left|\begin{array}{cccc}
c_1&c_2&\cdots&c_{n}\\
c_2&c_3&\cdots&c_{n+1}\\
\cdots&\cdots&\cdots&\cdots\\
c_{n}&c_{n+1}&\cdots&c_{2n-1}\end{array}\right|,
\end{split}
\end{align}
which correspond $A_1=\alpha_{0}, A_2=\alpha_2, \cdots,\ B_1=\alpha_{1}, B_2=\alpha_3,\ \cdots$, 
that is, $A_n=\alpha_{2n-2},\ B_n=\alpha_{2n-1}$ where $\alpha_m$'s are 
our Muir-Rogers determinants. 
Van Vleck (and Stieltjes) claims $A_n, B_n >0$ in the context of probability theory, the proof of which is omitted here. 
The coefficients of corresponding continued fraction are therefore given by
\begin{align}
&e_0=\alpha_0,\ e_1=\frac{\alpha_1}{\alpha_0},\ e_2=\frac{\alpha_2}{\alpha_1\alpha_0},\ 
e_3=\frac{\alpha_3\alpha_0}{\alpha_2\alpha_1},\ \cdots,
\nonumber \\
&e_n=\frac{\alpha_n\alpha_{n-3}}{\alpha_{n-1}\alpha_{n-2}}\quad(n\geq 3),
\label{VanVleckF}
\end{align}
which are our \eqref{MuirRogersF}. \\

\noindent
{\bf (Quotient-Difference algorithm)}\quad 
Now we wish to remark here that the right hand side of \eqref{VanVleckF} reminds us 
another theory, {\it i.e.} the quotient-difference (QD) algorithm by Rutishauser \cite{Rutishauser}, 
which is a method to find eigenvalues of a given matrix $A$. Let us give a short introduction of QD method. 
The eigenvalues of $N\times N$ matrix $A$, assumed to be symmetric $A^{\text{T}}=A$ for simplicity, 
are given by the zeros of characteristic polynomial $P(z)$, defined by
\begin{align}
P(z)=\text{det}\left(zE-A \right)=a_0z^N+a_1z^{N-1}+\cdots+a_N,
\end{align}
where $a_0=1$. Then we introduce a power series
\begin{align}
f(x)&=\frac{z^N}{P(z)}=\frac{1}{a_0+a_1x+\cdots+a_Nx^N}\qquad (x=z^{-1})
\nonumber \\
&=b_0+b_1x+b_2x^2+\cdots,
\label{b_s}
\end{align}
whose poles are the inverse of eigenvalues of $N\times N$ matrix $A$. 
Without loss of generality we can assume non-zero eigenvalues and $b_0=1$. 

In 1993 the present author found \cite{SogoQD} that the procedure of QD algorithm is 
equivalent to discrete dynamics of Toda molecule equation \cite{Hirota, Toda}. 
Let us suppose the matrix $A$ is given by a tri-diagonal matrix
\begin{align}
A=\left(\begin{array}{ccccc}
a_1&b_1&0&\cdots&0\\
b_1&a_2&b_2&\cdots&0\\
\cdots&\cdots&\cdots&\cdots&\cdots\\
0&\cdots&b_{N-2}&a_{N-1}&b_{N-1}\\
0&\cdots&0&b_{N-1}&a_N\end{array}\right)=A(0).
\label{matrixA}
\end{align}
Since any real symmetric matrix can be tri-diagonalized by Householder method, we can start from this initial condition.

Then \eqref{matrixA} is decomposed as $A=LR$ by
\begin{align}
\begin{split}
&L=\left(\begin{array}{ccccc}
1&0&0&\cdots&0\\
V_1&1&0&\cdots&0\\
\cdots&\cdots&\cdots&\cdots&\cdots\\
0&\cdots&V_{N-2}&1&0\\
0&\cdots&\cdots&V_{N-1}&1
\end{array}\right)=L(0),\\
&R=\left(\begin{array}{ccccc}
I_1&r_1&0&\cdots&0\\
0&I_2&r_2&\cdots&0\\
\cdots&\cdots&\cdots&\cdots&\cdots\\
0&\cdots&\cdots&I_{N-1}&r_{N-1}\\
0&\cdots&\cdots&0&I_N
\end{array}\right)=R(0),
\end{split}
\end{align}
which is a special case of so-called $LU$ decomposition.
The equality $A=LR$ implies in terms of matrix elements
\begin{align}
\begin{split}
&a_1=I_1,\ a_2=I_2+r_1V_1,\ \cdots,\ a_N=I_N+r_{N-1}V_{N-1},
\\
&b_1=r_1=V_1I_1,\ \cdots,\ b_{N-1}=r_{N-1}=V_{N-1}I_{N-1},
\end{split}
\end{align}
from which the initial values $I_1,\cdots,I_N$ and $V_1,\cdots,V_{N-1}$ are determined. 
The variables $(I_n, V_n)\ (n=1,2,\cdots, N)$ are pairs of canonical variables of discrete time Toda molecule system, 
where $V_N=0$ is assumed which is the reason of the name {\it molecule} instead of lattice.

Now we consider the following iteration, which is called $LR$ (left-right) algorithm, with $\ell=0,1,2,\cdots$,
\begin{align}
A(\ell+1)=L(\ell+1)R(\ell+1)=R(\ell)L(\ell),
\end{align}
which is written in matrix elements by
\begin{align}
\begin{split}
&r_n(\ell+1)=r_n(\ell)\equiv r_n,\\
&I_n(\ell+1)-I_n(\ell)=r_nV_n(\ell)-r_{n-1}V_{n-1}(\ell+1),\\
&I_n(\ell+1)V_n(\ell+1)=I_{n+1}(\ell)V_n(\ell).
\end{split}
\label{TodaMolecule}
\end{align}
This is the QD algorithm and also the Toda molecule equation with discrete time $\ell$. 
It should be noted that
\begin{align}
A(\ell+1)=L(\ell)^{-1}A(\ell)L(\ell),
\end{align}
which is orthogonal transformation, and keeps eigenvalues unchanged.  

In the QD algorithm it is convenient to introduce {\it tau function} $\tau_n(\ell)$, which is a determinant defined by
\begin{align}
\tau_n(\ell)=\left|\begin{array}{cccc}
b_\ell&b_{\ell+1}&\cdots&b_{\ell+n-1}\\
b_{\ell+1}&b_{\ell+2}&\cdots&b_{\ell+n}\\
\cdots&\cdots&\cdots&\cdots\\
b_{\ell+n-1}&b_{\ell+n}&\cdots&b_{\ell+2n-2}\end{array}\right|,
\label{TodaMoleculeS}
\end{align}
which is again Hankel type determinant with matrix elements given by coefficients $b_\ell$'s of \eqref{b_s}. 
In terms of $\tau_n(\ell)$, Toda molecule equation \eqref{TodaMolecule} are solved by 
\begin{align}
I_n(\ell)=\frac{\tau_{n}(\ell+1)\tau_{n-1}(\ell)}{\tau_{n}(\ell)\tau_{n-1}(\ell+1)},\quad
r_nV_n(\ell)=\frac{\tau_{n+1}(\ell)\tau_{n-1}(\ell+1)}{\tau_{n}(\ell)\tau_{n}(\ell+1)},
\end{align}
where only one condition remained is the following Hirota equation
\begin{align}
\tau_n(\ell+1)\tau_n(\ell-1)-\tau_n(\ell)^2=\tau_{n+1}(\ell-1)\tau_{n-1}(\ell+1),
\end{align}
and this is solved by \eqref{TodaMoleculeS} automatically.
 
After large enough iterations, $L(\infty)$ becomes unit matrix, and therefore we obtain
\begin{align}
A(\infty)=\left(\begin{array}{ccccc}
I_1(\infty)&r_1&\cdots&\cdots&0\\
0&I_2(\infty)&r_2&\cdots&0\\
\cdots&\cdots&\cdots&\cdots&\cdots\\
0&\cdots&\cdots&I_{N-1}(\infty)&r_{N-1}\\
0&\cdots&\cdots&0&I_N(\infty)\end{array}\right),
\end{align}
which implies $I_1(\infty), I_2(\infty),\cdots,I_N(\infty)$ are the eigenvalues to be determined. 
Physically speaking, this is the phenomenon that solitons with different characteristics $I_n$'s, 
which correspond to {\it velocities of solitons}, decouple after long time progress. 

The reason why $L(\infty)=E$ is as follows. Let $N$ poles of $f(x)$ are $x_1,x_2,\cdots,x_N$ with 
$|x_1|<|x_2|<\cdots<|x_N|$, we can show \cite{Rutishauser}
\begin{align}
\tau_n(\ell) \fallingdotseq \frac{C_n}{(x_1x_2\cdots x_n)^\ell},
\end{align}
from which we have
\begin{align}
V_n(\ell)\fallingdotseq\frac{C_{n+1}C_{n-1}}{C_n^2} x_n\left(\frac{x_n}{x_{n+1}}\right)^\ell\rightarrow 0
\quad\text{as}\quad \ell\rightarrow\infty.
\end{align}

Thus far we have noticed that QD algorithm is very similar to Muir-Rogers formula, except the former has an additional parameter 
$\ell$, which is {\it time} or iteration step. Such extension might be realized by an inclusion of $\ell$ into $\Phi(u)$ 
such as $\Phi_\ell(u)$, a time-dependent theory of probability distribution, or orthogonal polynomials \cite{SogoTDOP}. 
Unfortunately the author has not arrived at a decisive conclusion yet, so let us postpone it to future work.\\

\noindent
{\bf (Modular property)}\quad
On account of \eqref{Elkies}, there is the following fact for $a=7$. If we set
\begin{align}
\begin{split}
&x=-q^{4/7}\prod_{m=1}^\infty (1-q^{7m})(1-q^{7m-1})(1-q^{7m-6}),
\\
&y=+q^{2/7}\prod_{m=1}^\infty (1-q^{7m})(1-q^{7m-2})(1-q^{7m-5}),
\\
&z=+q^{1/7}\prod_{m=1}^\infty (1-q^{7m})(1-q^{7m-3})(1-q^{7m-4}),
\end{split}
\label{ElkiesF}
\end{align}
they satisfy the equation
\begin{align}
x^3y+y^3z+z^3x=0,
\label{Klein}
\end{align}
which is called Klein's quartic curve in CP$^2$ (two-dimensional complex projective space, see p.84 of \cite{Elkies}). 
It might be possible to extend \eqref{Klein} to Kummer surface \cite{Hudson} by introducing additional parameters in \eqref{ElkiesF}. 
There might be other models of curve in higher dimensions by choosing larger values of $a$.

%
\end{document}